\documentstyle[12pt]{amsart}
\newtheorem{thm}{\bf Theorem}[section]
\newtheorem{cor}[thm]{\bf Corollary} 
\newtheorem{lem}[thm]{\bf Lemma} 
\newtheorem{prop}[thm]{\bf Proposition} 
\newtheorem{definition}[thm]{\bf Definition} 
 
\newtheorem{rem}[thm]{\bf Remark} 
 
\newtheorem{conj}[thm]{\bf Conjecture}

\def\proof{\noindent {\sl Proof:}\;  }
\def\t{\noindent} 

\def \jeden {1\hskip-3.5pt1}

\def\codim{ \mbox{\rm codim}\, }
\def\Dual{\mbox{\rm Dual}\, }

\def\K{{\cal K}}
\def\F{{\cal F}}

\def\Z{{\Bbb Z}}
\def\Q{{\Bbb Q}}

\def\C{{\Bbb C}}
\def\P{{\Bbb P}}

\newcommand{\maprU}[1]{%  %%%%%%%%%% ___-____> %%%%%
\smash{\mathop{%
\hbox to 1cm{\rightarrowfill}}\limits^{#1}}}
\newcommand{\maprD}[1]{%  %%%%%%%%%% ---_----> %%%%%
\smash{\mathop{%
\hbox to 1cm{\rightarrowfill}}\limits_{#1}}}
\newcommand{\maplU}[1]{%  %%%%%%%%%% <___-____ %%%%%
\smash{\mathop{%
\hbox to 1cm{\leftarrowfill}}\limits^{#1}}}
\newcommand{\maplD}[1]{%  %%%%%%%%%% <---_---- %%%%%
\smash{\mathop{%
\hbox to 1cm{\leftarrowfill}}\limits_{#1}}}

\newcommand{\mapdL}[1]{%  %%%%%%%%%  V*VV %%%%%%%
\downarrow\llap{$\vcenter{\hbox{$\scriptstyle#1\,\,\,$}}$}}

\begin{document}
%\sf 
%
%
\title[Equivariant Chern-Schwartz-MacPherson class]
{Equivariant Chern classes of singular algebraic varieties
with group actions}
\author[T.~Ohmoto]{Toru Ohmoto}
\thanks{Partially supported by Grant-in-Aid for Encouragements of
Young Scientists (No.15740042), the Japanese Ministry of Education, Science
and Culture. }
\address[T.~Ohmoto]{Department of Mathematics, 
Faculty of Science,  Hokkaido University,
Sapporo 060-0810, Japan}
\email{ohmoto@@math.sci.hokudai.ac.jp}
\keywords
{Equivariant Chern class, Chern-Schwartz-MacPherson class, Equivariant
homology, Classifying space, Thom polynomial  }
%
%
%\date{\today}
%
\maketitle
\begin{abstract}
We define the equivariant Chern-Schwartz-MacPherson class of 
a possibly singular algebraic $G$-variety over the base field $\C$, 
or  more generally over  a field of characteristic $0$. 
In fact, we construct a  natural transformation $C^G_*$ 
from the $G$-equivariant constructible function functor $\F^G$  to 
the $G$-equivariant homology functor $H^G_*$ or $A^G_*$  
(in the sense of Totaro-Edidin-Graham). 
This $C^G_*$ may be regarded 
as MacPherson's transformation for (certain) quotient stacks. 
We discuss on other type Chern classes  and applications. 
The Verdier-Riemann-Roch formula takes a key role throughout. 
\end{abstract}

%\setlength{\baselineskip}{16pt}

%%%%%%%%%%%%%%%%%%%%%%%%%%%%%%%%%%%%%%%%%%%
\section{Introduction}
For a possibly  singular complex algebraic  variety  $X$ there are 
several kinds of  ``Chern classes" of $X$ available. 
These ``Chern classes" of $X$  live in appropriate {\it homology groups} of
$X$, 
which satisfy ``the normalization property" that if $X$ is non-singular, then 
it coincides with the Poincar\'e dual to the ordinary Chern class 
of the tangent bundle $TX$.   

 {\it The Chern-Schwartz-MacPherson class}  is one of them. 
R.~MacPherson \cite{Mac} constructed the class 
 to solve the so-called Grothendieck-Deligne conjecture: 
Actually  he proved the existence of 
a unique natural transformation  $C_*:\F(X) \to H_{2*}(X ;\Z)$ 
from  the abelian group $\F(X)$ of {\it constructible functions} over $X$ 
to the homology group (of even dimension)
so that if $X$ is nonsingular, then $C_*(\jeden_X) =c(TX)\frown [X]$ 
where  $\jeden_X$ is the characteristic function, $\jeden_X(x)=1$ ($x \in X$).  
Independently M.~H.~Schwartz \cite{Sch}  
had introduced obstruction classes (defined in a local cohomology)  
for the extension of stratified radial vector frames over $X$, and 
it is shown (\cite{BS}) that  both classes coincide, 
so $C_*(\jeden_X)$ is often denoted by $C^{SM}(X)$.  
In a purely algebraic context, 
MacPherson's transformation is also formulated as 
$C_*:\F(X) \to A_*(X)$, 
the value being in the Chow group of cycles modulo rational equivalence, 
for embeddable  schemes (separated and of finite type) 
over arbitrary base field $k$ of characteristic $0$. 
That was done by G.~Kennedy \cite{Ken} 
using the groups of {\it Lagrangian cycles} (cf.  \cite{Gon},
\cite{Sab}), 
which is isomorphic to $\F(X)$ in a certain way. 
In the complex analytic context, 
MacPherson's  theory  is also verified: for instance, 
the crucial step in \cite{Mac}, {\it the graph construction}, 
is proved in \cite{Kw1} in the analytic setting. 
Besides, the Lagrangian cycle approach in the complex differential geometry
\cite{Fu}  
and Schwartz's approach within the Chern-Weil theory \cite{BLSS}  
have been also achieved.  

In this paper we  think of  a $G$-version 
of the Chern-Schwartz-MacPherson class
 for algebraic $G$-varieties $X$.   
Our main aim is to focus 
the elementary  (or formal)  construction of the equivariant version of $C_*$ 
as well $G$-versions of $\F(X)$ and $H_*(X)$ (or $A_*(X)$). 
So for the simplicity  
we  discuss  basically in the complex context like as the original \cite{Mac}: 
Then we use the singular cohomology and 
 the Borel-Moore homology, simply denoted by $H_*(X)$, 
of the underlying {\it analytic} space  
(denoted by the same letter $X$ for short). 
However, after suitable changes, 
 the reader can  read them as in the algebraic context (\cite{Ken}) 
with the use of  (operational) Chow rings and Chow groups: 
Then a  {\it  scheme} is assumed to be  separated  and of finite type 
over  $k$ of characteristic $0$, 
and a {\it variety}  is an irreducible and reduced such scheme. 

As known, for a topological group $G$, 
Borel's equivariant cohomology of a $G$-space $X$ is defined by 
$$H_G^*(X)=H^*(X \times_G EG),$$
where $EG \to BG$ is the universal principal bundle over the classifying
space of $G$. 
As a counterpart in algebraic geometry, 
for reductive linear algebraic group $G$, 
the {\it $G$-equivariant homology group} $H^G_*(X)$  
($A^G_*(X)$) of a $G$-variety $X$ 
is defined in Edidin-Graham \cite{EG} 
using the algebraic approximation of $BG$ given by Totaro \cite{Totaro}.  
From the same viewpoint,  we introduce the abelian group 
$\F^G(X)$ of {\it $G$-equivariant constructible functions} over $X$ 
(that is, roughly,  constructible functions over  $X \times_G EG$ 
whose supports have finite codimension). In particular the group 
$\F^G_{inv}(X)$ of {\it $G$-invariant} constructible functions over $X$ 
becomes a subgroup of $\F^G(X)$  by  a natural identification. 
Both of $\F^G$ and $H^G_*$  become {\it covariant} functors for the category of 
$G$-varieties and proper $G$-morphisms (see subsections \ref{eqhomology}
and \ref{eqconstf}). 

From now on 
we assume that 
a $G$-variety (scheme) $X$  has   a closed equivariant embedding into 
some $G$-nonsingular varieties, 
and when we emphasize it, we  say such $X$ is {\it $G$-embeddable} for short.
We show the following theorem for $G$-embeddable varieties: 

\begin{thm}\label{C}
Let $G$ be a complex reductive linear algebraic group. 
For the category of complex algebraic $G$-varieties $X$ and proper
$G$-morphisms, 
there is a  natural transformation of covariant functors 
$$C_*^G: \F^G(X) \to H^G_*(X)$$
such that if $X$ is non-singular, then $C_*^G(\jeden_X)=c^G(TX) \frown [X]_G$ 
where $c^G(TX)$ is $G$-equivariant total Chern class of the tangent bundle
of $X$. 
The natural transformation $C^G_*$ is unique in a certain sense. 
\end{thm}

To be precise, we mean by  {\it a natural transformation} that 
$C^G_*$ satisfies that $(i)$:
$C^G_*(\alpha+\beta)=C^G_*(\alpha)+C^G_*(\beta)$ and 
$(ii)$: $f^G_*C^G_*=C^G_*f^G_*$ for any proper $G$-morphism $f: X \to Y$.  
The precise statement of  the ``uniqueness" of $C^G_*$ is seen in the
subsection 3.2 (b).  

\begin{rem} \label{remC} 
Theorem \ref{C} is also true over the base field $k$ of characteristic $0$, 
(at least) for quasi-projective schemes $X$ with linearlized $G$-action; 
then we have a natural transformation $C_*^G: \F^G(X) \to A^G_*(X)$ 
which satisfies the normalization property 
(see the subsection \ref{mixed}  and the proof of Theorem \ref{C} given in
\S 3).  
This $C^G_*$ is naturally regarded 
as the extension of MacPherson's Chern class theory to the category of
quotient stacks,  
$C_*:\F([X/G]) \to A_*([X/G])$ (Theorem \ref{Cq}). 
\end{rem} 

\begin{definition}
The $G$-equivariant Chern-Schwartz-MacPherson class of a $G$-variety $X$ 
is defined by   $C^{SM}_G(X):=C_*^G(\jeden_X)$. 
\end{definition}
Besides of $\jeden_X$ and $C^{SM}_G(X)$, 
we can take other kinds of ``canonical  constructible functions" over $X$ and 
the corresponding ``canonical Chern-SM classes" of $X$. That will be
discussed in \S 6.

The rest of this paper is organized as follows: 

In \S 2  we will  review some  basic materials from \cite{Totaro} and \cite{EG} 
but in a slightly different form. 
Groups $\F^G(X)$ and $H^G_*(X)$ ($A^G_*(X)$)  are defined 
to be the inductive limit of abelian groups 
via very simple ``Radon transforms" (labeled by $(*_F)$ and $(*_H)$). 
In \S 3  our $C^G_*$ is given as the limit of 
``MacPherson's  transformation for topological Radon transforms" studied in
\cite{EOY}. 
Then Theorem \ref{C} automatically follows.  
We remark that this construction is very related to the ``proconstruction"
of $C_*$ 
for {\it provarieties} (projective limits of varieties) given by Yokura
\cite{Y2} 
(Remark \ref{RemarkLim}).  
The last subsection \ref{Chernq} of \S 3 is devoted to 
the interpretation of $C^G_*$  in terms of quotient stacks. 

In \S 4 we note some useful properties of our $C^G_*$, for instance, 
the equivariant versions of {\it Verdier-Riemann-Roch formula} 
(written by VRR formula for short) 
 for smooth morphisms (\cite{FM}, \cite{Y1},  
and also \cite{Schurmann} for local complete intersection morphisms). 
That is a Riemann-Roch type theorem saying 
the compatibility of the transformation $C_*$ with 
certain {\it pullbacks}  (i.e., contravariant operation) 
of constructible functions and homologies. 
In fact, the simplest VRR formula is involved in our construction of
$C^G_*$ itself  
(the square (2) in the proof of Lemma \ref{com12}). 

 {\it The equivariant Chern-Mather class} $C^M_G(X)$ is introduced in \S 5. 
As known, 
the Chern-Mather class $C^M(X)$ is 
a key factor in the construction of  (ordinary) MacPherson's class,  
which is roughly  
the Chern class of limiting tangent spaces of the regular part of $X$. 
In fact $C^{SM}(X)$ is expressed by $C^M(X)$ plus 
a certain linear combination of $C^M(W)$'s of subvarieties $W$  in the
singular locus of $X$. 
We show the equivariant version of such relations. 
Then some simple properties of $C^G_*$ become clearer: 
for instance,  the restriction of  $C^G_*$ 
to a fibre  of the universal principal bundle $X\times_G EG\to BG$ 
recovers the ordinary MacPherson transformation $C_*$.  

In \S 6 and \S 7, we discuss on  apparently a bit different  two kinds of applications, 
which generalize 
 {\it orbifold Euler characteristics} (cf.   \cite{HH}, \cite{BF}) 
and {\it Thom polynomials} (cf.  \cite{Thom}, \cite{Kaz1}, \cite{FR}), respectively.  
As to the former topic, {\it the canonical quotient Chern classes}  are introduced, 
which reflect some commutator structure of the group action.  
In \cite{O2} we will apply this theory to typical examples such as symmetric products, 
and obtain generating functions of the quotient Chern classes whose constant terms 
provide well-known generating functions of (orbifold) Euler characteristics. 
As to the latter topic,  a Thom polynomial is roughly saying 
the $G$-Poincar\'e dual to an invariant subvariety of a $G$-nonsingular variety. 
As a simple generalization, we study 
the $G$-Poincar\'e dual to the ``Segre-version" of our equivariant Chern class. 
Our Theorem \ref{segretp} is motivated by 
the formula of Parusi\'nski-Pragacz \cite{PP} (Theorem 2.1) 
for degeneracy loci of generic vector bundle morphisms. 

Consequently,  it can be viewed that 
these two applications deal with a unified  ``Chern class version" 
of  {\it Euler characteristics} and {\it fundamental classes}  
arising in  some ``$G$-classification theory".

The author would like to thank  especially Shoji Yokura for discussions and 
his comments on the first draft of this paper.

%%%%%%%%%%%%%%%%%%%%%%%%%%%%%%%%%%%%%%%%%%%

\section{Classifying space and the Borel construction}
In 2.1 -- 2.5  we pick up some definitions and properties from
\cite{Totaro} and \cite{EG}
Note again that $H_*$ can be replaced by the Chow group $A_*$ in the
algebraic context. 
In 2.6 the equivariant constructible function is defined. 

\subsection{Totaro's construction of $BG$}
Let $G$ be a  complex reductive linear algebraic group of dimension $g$. 
Take an $l$-dimensional representation $V$ of $G$  
with a $G$-invariant Zariski closed subset $S$  in $V$ so that $G$ acts on
$U:=V-S$ freely. 
It is possible to take  $V$ and $S$  so that 
the quotient $U\to U/G$ becomes 
an algebraic principal $G$-bundle over a quasi-projective variety, 
and that the codimension of $S$ is  sufficiently high.   
Actually this is achieved by a similar construction of Grassmanian varieties 
(Remark 1.4 of \cite{Totaro}). 
Let $I(G)$  be the collection of  Zariski open sets $U=V-S$ 
where $V$ is a representation and $S$ is a closed subset  of $V$ with
properties just as mentioned. 
We put a partial order on $I(G)$: we say $U(=V-S) < U'(=V'-S')$  
if $\codim_V\, S  <  \codim_{V'}\, S'$  
and there is an $G$-equivariant linear  inclusion $ V \to V'$ 
sending $U$ into $U'$.  Then $(I(G), <)$ is a directed set. 
All quotients $U \to U/G$ with induced maps by inclusions  form  an
inductive system, that is 
the algebraic approximation of the universal principal bundle $EG \to BG$
(\cite{Totaro}, \cite{EG}). 

An algebraic construction of classifying maps 
(for principal bundles over quasi-projective varieties)  is given in 
 Lemma 1.6 in \cite{Totaro}, that will be 
used in the last section. 

\subsection{Mixed quotients} \label{mixed} 
Let $X$ be a $G$-variety.  
For any $U\in I(G)$, 
the diagonal action of $G$ on $X \times U$, which is always a free action,  
gives a principal bundle $X \times  U \to X\times_G  U=(X\times  U)/G$, 
and thus  
the equivariant projection $X\times  U \to U$ serves  
the fibre bundle $X\times_G  U \to U/G$  with fibre $X$. 
Roughly saying, 
the universal fibre bundle  $X \times_G EG \to BG$ is approximated by 
those mixed quotients $X\times_G  U \to U/G$ for $U \in I(G)$
(Edidin-Graham \cite{EG}).  

We attention to the fact that 
in general the mixed quotients $ X\times_G  U$ exists 
as {\it algebraic spaces} in the sense of Artin, not as schemes
(Proposition 22 \cite{EG}).  
To avoid this technicality, 
we may think of the following cases: 
In the complex case $k=\C$, a separated algebraic space  of finite type  admits 
the corresponding analytic space (Corollary 1.6 in \cite{Artin}) 
(besides, in our convention $X$ is assumed to be  separated, of finite type
and ($G$-)embeddable, 
hence $X\times_G U$ is also as an algebraic space, thus as an analytic space). 
In the algebraic context, 
we assume that $X$ is a quasi-projective scheme with a linearlized $G$-action. 
Then the  mixed quotient $X\times_G U$ 
exists  as a  quasi-projective scheme (Proposition 23 in  \cite{EG}). 
Note that this quasi-projective hypothesis covers rather many interesting
cases. 

In fact,  we will later apply (ordinary) transformation $C_*$ 
to those mixed quotients in the complex case (by appealing to the
transcendental method)  
and also in the quasi-projective case over $k$ of characteristic $0$. 
Presumably Kennedy's formulation (for embeddable schemes) would be extendable 
into the context of (embeddable) algebraic spaces, 
then this kind restriction mentioned above would not be needed.

\subsection{$G$-equivariant cohomology} 
For any pair $U$ and $U'$ ($\in I(G)$)  so that $U<U'$, 
we let $\iota_{U,U'}: X \times_G U \to X \times_G U'$ denote the natural
inclusion 
and $r_{U',U}:=\iota_{U,U'}^*: H^*(X\times_G  U') \to H^*(X\times_G  U)$ 
the induced homomorphism. 
Then we have a projective system $\{H^*(X\times_G  U), r_{U,U'}\}$ and 
{\it the $i$-th equivariant cohomology} of $X$ is given as 
$$H_G^i(X)=\lim_{\underset{I(G)}{\longleftarrow}} H^i(X\times_G  U).$$ 
The formal sum is denoted by $H_G^*(X)=\prod H_G^i(X) =\lim_{\leftarrow}
H^*(X\times_G  U)$. 
We also denote by $r_U: H_G^*(X) \to H^*(X\times_G  U)$ the canonical
projection for $U$.  
Note that in the algebraic context, 
the  cohomology groups should be replaced 
by {\it operational Chow groups} (\cite{F}, \cite{EG}). 

Let $\xi$ be a $G$-equivariant vector bundle $E \to X$
 (i.e., $E$, $X$ are $G$-varieties, the projection is $G$-equivariant), 
then $\xi$ induces a vector bundle 
$E \times_G U \to X \times_G  U$, denoted by $\xi_{U}$. 
The projective limit of Chern classes $c(\xi_{U})$ gives 
{\it the $G$-equivariant Chern class} of  $\xi$, which is denoted by
$c^G(\xi)  \in H_G^*(X)$. 
In particular, when $X=\{pt\}$, an equivariant vector bundle is $V \to
\{pt\}$  being $V$
a representation. The Chern class is denoted by $c^G(V)\in H_G^*(pt)=H^*(BG)$.  
Besides, the pullback of $c^G(V)$ via the trivial equivariant morphism 
$X \to  \{pt\}$ is denoted by the same notation: $c^G(V) \in H^*_G(X)$. 

Given a $G$-morphism $f: X \to Y$, the pullback $f_G^*: H_G^i(Y) \to
H_G^i(X)$ is defined 
in a natural way: $f_G^*(\{\alpha_U\}) :=\{(f\times_G id)^*\alpha_U\}$. 

\subsection{$G$-equivariant homology} \label{eqhomology} 
We repeat Edidin-Graham's definition of the equivariant homology (Chow group) 
(Proposition 1 in \cite{EG})   
but in a suitable form for the later use 
(again we describe it in the complex context but it works over any case).

At first we define a sub-order $<_*$ on $I(G)$: 
For any two $U(=V-S)$ and $U'(=V'-S')$, 
we say that 
$U<_* U'$  if  there is a representation $V_1$ so that 
$V\oplus V_1=V'$ and $U\oplus V_1 \subset U'$.  
Note that if $U_1<U_2$, then there is $U'$ so that  
$U_1<_* U'$ and $U_2<_* U'$ (e.g., $U'=V_1\oplus V_2-S_1\oplus S_2$). 

Let $X$ be a complex $G$-variety of $\dim X=n$ (equidimensional). 
To each $U=V-S$ with $\dim V=l$ and $\codim  S = s$, we assign a {\it
truncated} homology  
$$H_{trunc}(X\times_G U):=\underset{2(n-s)<i \le 2n}
{\oplus}H_{i+2(l-g)}(X\times_G U).$$ 
Note that the range of dimension in the direct sum depends on $U$ (precisely, 
the dimensions of $V$ and $S$). 
This notation is convenient for us  because 
 we shall later think of  {\it total} homology classes ({\it total Chern
classes})  
rather than a distinguished $i$-th homology class.

For each pair $U<_* U'$ ($V'=V\oplus V_1$), the diagram 
$U \leftarrow U\oplus V_1 \rightarrow U'$ 
of projection and injection induces  
$$
\begin{array}{rcl}
p=p_{U,U\oplus V_1}\quad & X\times_G (U\oplus V_1)&\quad \iota_{U\oplus
V_1, U'}=\iota \\
\swarrow& & \searrow\\
X\times_G U& & X\times_G U'
\end{array}
\eqno{(*)}
$$
This diagram $(*)$ induces the following isomorphisms for $2(n-s)< i\le 2n$ 
(we denote $(i^*)^{-1}$ by $\iota_*$, in abusing the notation): 
$$
\begin{array}{c}
\qquad \qquad 
p^* \quad  H_{i+2(l+k-g)}(X\times_G (U\oplus V_1)) \quad\iota_*:=(\iota^*)^{-1}
\\
\simeq  \, \nearrow  
\qquad \qquad \quad \qquad \qquad \searrow  \, \simeq 
\\
H_{i+2(l-g)}(X\times_G U) \qquad \qquad \qquad  
H_{i+2(l+k-g)}(X\times_G U')
\end{array}
\eqno{(*_H)}$$

This is because $\iota$ is an open embedding 
(so the pullback $\iota^*$ is defined) and its complement 
$X\times_G U'-X\times_G (U\oplus V_1)$  has 
the (complex) codimension $\ge s$  (hence $\iota^*$ is isomorphic). 
Also $p=p_{U,U\oplus V_1}$ is the projection of a vector bundle 
so it induces an isomorphism $p^*$. 
The composition of isomorphisms in $(*_H)$ define  a graded  homomorphism of 
truncated homology groups, whose degrees are shifted by $k=\dim U'-\dim U$,
 denoted by 
$$\varphi_{U,U'}: 
\underset{2(n-s)<i } {\oplus}H_{i+2(l-g)}(X\times_G U) \to 
\underset{2(n-s')<i} {\oplus}H_{i+2(l+k-g)}(X\times_G U').$$
This makes an inductive system with respect to the directed set $(I(G),<_*)$, 
moreover with the original order $<$ 
\footnote{
For any pair $U_1<U_2$ with respect to the original order,  take $U'$ so that  
$U_1, U_2<_* U'$ (e.g., $U'=U_1\oplus U_2$). 
Then  both $\varphi_{U_1,U'}$ and $\varphi_{U_2,U'}$ are isomorphic at least in 
the range  $i>2n-2s_1$ ($s_1=\codim S_1$), and hence 
we can define a canonical injective homomorphism $\varphi_{U_1,U_2}$ 
($:=(\varphi_{U_2,U'})^{-1}\circ \varphi_{U_1,U'}$) 
from the truncated homologies of $U_1$ to the one of $U_2$, that is 
so-called the double filtration argument in \cite{Totaro}, \cite{EG}. 
Therefore we don't take care of underlying orders, so we  use $<_*$ 
basically. }. 
We define {\it the $i$-th equivariant homology group}  to be  (the limit of) 
the shifted-dimensional component of the truncated homology 
$$H^G_i(X) = H_{i+2(\dim U-g)}(X\times_G  U) \quad 
\mbox{(for $U$ with $\codim S$ high enough)}.$$
Thus $H_i^G(X)$ is trivial for $i>2n$ and  possibly nontrivial for any
negative $i$. 
The direct sum is denoted by 
$$H_*^G(X)=\oplus  H^G_i(X)=\lim_{\longrightarrow} H_{trunc}(X\times_G  U).$$
For each $U$, the identification map is denoted 
by $\varphi_U: H_{trunc}(X \times_G U) \to H_*^G(X)$.  

Given a proper $G$-morphism $f: X \to Y$ between $G$-varieties, we have 
an induced homomorphism  $f^G_*: H_*^G(X) \to H_*^G(Y)$ defined by 
$f^G_*(\varphi_U(c)) :=\varphi_U((f\times_G id)_*(c))$ 
as the limit of $(f\times_G id)_*: H_{trunc}(X\times_G U) \to
H_{trunc}(Y\times_G U)$. 
Any other expected functorial properties  are also  satisfied, see \cite{EG}.  

As an example we illustrate a most simplest case: $G=GL(1)$, $X=\{pt\}$ and 
a sequence in $I(G)$  that is $\{U_m=\C^{m+1}-\{0\}\}$ with the action of
all weights $-1$: 
$$
\begin{array}{rcl}
&(U_m\times \C)/GL(1)& =\P^{m+1}-\{pt\}\\
\cdots \;\;\stackrel{p}{\swarrow}\;\;\searrow\;\;\stackrel{p}{\swarrow} 
&&\searrow \;\; \cdots \\
\P^1\subset \cdots \subset \;\;  \P^m\;\; 
&\subset&\;\; \P^{m+1} \;\;\subset \cdots \subset \P^\infty \sim BGL(1)
\end{array}
$$
Note that $n=\dim X=0$, $g=\dim G=1$ and 
$$H_{trunc}(X\times_G U_m)
=\underset{-2(m+1)<i}{\oplus} H_{i+2m}(\P^m) \simeq \Z^{m+1}.$$
The map $\varphi_{U_m, U_l}:\Z^{m+1} \to \Z^{l+1}$ is a natural inclusion 
and hence $H^{GL(1)}_{i}(pt)=\Z$ for nonpositive even number $i$, and 
trivial otherwise.

Roughly saying,  the direct sum of fibres of $p$'s over a point approximates  
``the tangent space $T_{\P^\infty}$". 
 In fact there will appears somewhat  ``(inverse) Chern class factor of
$T_{\P^\infty}$" 
in our definition of  $C^G_*$ given in \S 3 (Remark \ref{RemarkLim}).

\subsection{$G$-fundamental class and Poincar\'e duality} 
For any $U$,  the fundamental cycle $[X\times_G  U]$  tends to 
a unique element of $H_{2n}^G(X)$, denoted by $[X]_G$.
This is called {\it the $G$-equivariant fundamental  class} of $X$. 
Note that  we can identify 
$H_{2j}^G(X) =H_{2(j+l)}^G(X\times V)$ for any representation $V$ 
($\dim V =l$) 
through the pullback isomorphisms induced by 
$(X\times V)\times_G U \to X\times_G U$. 
Let $W$ be an $(j+l)$-dimensional $G$-invariant reduced closed subscheme of
$X \times V$ and $i: W \to X \times V$  the $G$-inclusion. 
Then $W$ represents an equivariant homology class of $X$ as 
$i^G_*([W]_G) \in H_{2(j+l)}^G(X\times V)= H_{2j}^G(X)$. 
 {\it An $j$-dimensional $G$-equivariant algebraic cycle class} means 
a finite sum $\sum_k a_k \, (i_k)_*([W_k]_G) \in H_{2j}^G(X)$, 
where each $W_k$ is  an $(j+l)$-dimensional $G$-invariant subvariety of
some $X \times V$.  
Here  $j$ is possibly negative ($-l \le j \le n$). 
Of course,  a $G$-invariant cycle of $X$ represents 
an equivariant cycle class of nonnegative dimension.

There is a well-defined homomorphism 
$$\frown [X]_G: H_G^{2n-i}(X) \to H_{i}^G(X), \quad 
a \mapsto \varphi_U(r_U(a) \frown [X\times_G U]).$$
If $X$ is nonsingular, this is isomorphic for each $i$, called  
{\it the $G$-equivariant Poincar\'e dual}. 
In particular, when $X$ is a point, $H^G_{-k}(pt) \simeq H_G^k(pt) =H^k(BG)$. 
We denote by ${\rm Dual}_G$  the inverse of the map $\frown [X]_G$ (for
each $i$). 
The composite map  $r_U\circ {\rm Dual}_G \circ \varphi_U$ 
coincides with the ordinary Poincar\'e dual of $X\times_G U$ on the
truncated homology.

\subsection{$G$-equivariant constructible functions} \label{eqconstf} 
A constructible function over a complex algebraic variety $X$ is 
an integer valued function $\alpha: X \to \Z$ 
which has a finite partition of $X$ into constructible subsets so that 
the value of $\alpha$ is constant over each of the constructible sets. 
We let $\F(X)$ denote the Abelian group consisting of all constructible
functions over $X$. 
Any constructible function $\alpha \in \F(X)$ is represented by 
$\alpha = \sum_{i=1}^k \, a_i\jeden_{W_i}$ for some integers $a_i$ and 
subvarieties $W_i$ of $X$.  Here $\jeden_W$ denotes  the function taking values 
$1$ for $x \in W$ and $0$ otherwise. 
For any  proper morphism $f: X \to Y$, 
we define the pushforward $f_*:\F(X) \to \F(Y)$  by 
$$f(\alpha)(y):=  \sum_{i=1}^k \; a_i \, \chi(f^{-1}(y) \cap W_i)  
\quad (\alpha =\sum_{i=1}^k \, a_i \jeden_{W_i}, y \in Y), $$ 
where $\chi$ is the topological Euler characteristic 
with respect to the Borel-Moore homology groups. 
Note that $(f\circ g)_*=f_*\circ g_*$. 
Even if $f$ is not proper,  the sum in the right hand side may be finite for  
any $\alpha$ and $y$, and in that case we also denote the map by  $f_*$. 
For a constructible function $\alpha \in \F(X)$, 
we define {\it the integral of $\alpha$ over $X$} 
(or say, the Euler characteristics of $\alpha$)  to be the value
$f_*(\alpha) \in \F(pt)=\Z$ 
with $f:X\to \{pt\}$, that is 
$$\int_X \alpha := f_*(\alpha)=\sum_{i=1}^k \, a_i \,\chi(W_i) \; \in \; \Z. $$
For any morphism $f: X \to Y$, 
 the pullback $f^*:\F(Y) \to \F(X)$ is defined to be 
$f^*(\beta):=\beta\circ f$.  

\begin{rem}\label{chi}
In the case of the base field $k$ of characteristic $0$,  
the above definition of pushforward should be appropriately changed 
in terms of Lagrangian cycles, see \cite{Ken}. In abusing the notation,  
we may use the letter $\chi(X):=\int_X \jeden_X$ in this context too. 
\end{rem}

Now let $X$ be a variety with a $G$-action. 
The subgroup of $\F(X)$ consisting of 
$G$-invariant constructible functions is denoted by 
$$\F^G_{inv}(X) 
:= \{\; \alpha \in \F(X)\; | \; \alpha(g(x))=\alpha(x), \, (x \in X, g \in
G) \; \}.$$
For any $U<_* U'$  ($V'=V\oplus V_1$, $U=V-S$, $U'=V'-S'$),  
let $p:V' \to V$ be the projection to the first factor, 
then it induces a pullback homomorphism 
$$\phi_{U,U'}:=p^*:\F^G_{inv}(X\times  V) \to \F^G_{inv}(X\times  V'), \qquad 
\alpha \mapsto \alpha \circ (id \times p)$$
(we sometimes denote it by $\phi_{V,V'}$). 
Then, $\{\F^G_{inv}(X\times  V), \phi_{U,U'}\}$  makes an inductive system, so 
we define 
$$\F^G(X):=\lim_{\underset{I(G)}{\longrightarrow}}\F^G_{inv}(X\times V). $$ 
An element of this limit group is called  {\it a $G$-equivariant constructible
functions associated to $X$}. 
The limit map is denoted by  $\phi_U: \F^G_{inv}(X\times V) \to \F^G(X)$
(sometimes by $\phi_V$). 

In an obvious way, any $G$-invariant function over $X$ is lifted to 
an invariant function over $X \times V$, and hence 
there is a canonical inclusion, denoted by $\phi_0$, 
$$\F^G_{inv}(X) \subset\F^G(X), \quad \alpha \mapsto
\phi_0(\alpha)=\phi_V(\alpha \times \jeden_V).$$  
Note that if $X$ is a point, then $\F^G_{inv}(pt)=\F(pt)\simeq \Z$  
(consisting of constant functions) 
but $\F^G(pt)$ contains a lot of other invariant functions over
representations $V$'s. 

For a proper $G$-morphism $X \to Y$, we define 
{\it the equivariant pushforward} homomorphism  
$$f^G_*: \F^G(X) \to \F^G(Y), \quad 
f^G_*(\phi_U\alpha_U):=\phi_U((f\times id)_*(\alpha_U)), $$ 
that is the limit map of 
$(f\times id)_*: \F^G_{inv}(X\times V) \to \F^G_{inv}(Y\times V)$. 
It is easily checked that $(f\circ g)^G_*=f^G_*\circ g^G_*$ for proper
$G$-morphisms. 
Similarly {\it the equivariant pullback}  $f_G^*: \F^G(Y) \to \F^G(X)$ is
also defined.

For any $\alpha \in \F^G(X)$, 
we define {\it the $G$-integral} of $\alpha$ to be  $f^G_*(\alpha)  \in
\F^G(pt)$ 
by the pointed map $f:X \to \{pt\}$. 
So  a $G$-integral  is an equivariant constructible function over a point,
not constant in general.  
In particular, if $\alpha\in \F^G_{inv}(X)$, more precisely, 
$\alpha=\phi_V (\alpha_0 \times \jeden_V)$ for some $\alpha_0 \in
\F^G_{inv}(X)$,  
then the $G$-integral of $\alpha$  is the constant $\int_X \alpha_0$:  
\begin{eqnarray*}
f^G_*(\alpha) &=& f^G_*(\phi_V (\alpha_0 \times \jeden_V)) 
=\phi_V ((f\times id)_*(\alpha_0 \times \jeden_V)) \\
&=&\phi_V (f_*(\alpha_0) \times \jeden_V) = f_*(\alpha_0).
\end{eqnarray*}

\begin{rem} \label{Remarkconst}
Instead of $\phi_{U,U'}=p^*$, 
we think of  $\tilde{\phi}_{U,U'}$, the restriction of $\phi_{U,U'}$ to
$\F^G_{inv}(X\times U)$, 
which is regarded as the composed map 
$$
\begin{array}{rcl}
p^*=(p_{U,U\oplus V_1})^*\;\; &
\F^G_{inv}(X\times (U\oplus V_1))& \;\; (\iota_{U\oplus V_1, U'})_*=\iota_*\\
\nearrow& & \searrow  \\
\F^G_{inv}(X\times U)& & \F^G_{inv}(X\times U')
\end{array}
\eqno{(*_F)}
$$
($\iota$ is not proper but an open embedding, so $\iota_*$ is defined).   
Note that $\{\F^G_{inv}(X\times U), \tilde{\phi}_{U,U'}\}$ 
also makes an inductive system, 
but the limit group differs slightly from the above $\F^G(X)$.  
Later we will consider essentially this smaller limit group. The point is that 
 $\tilde{\phi}_{U,U'}$ for constructible functions and $\varphi_{U,U'}$ for
homology 
can be read off as  {\it topological Radon transforms}, 
that is, roughly saying, pulling back and then pushing forward. 

In \cite{EOY}   the authors studied 
the category  whose objects are  (nonsingular) varieties and 
``morphisms" of $X$ to $Y$ are diagrams 
$X \stackrel{p}{\leftarrow} M \stackrel{q}{\rightarrow} Y$ 
with  $p$ being a smooth morphism (the case of singular varieties is
supported by \cite{Y1}). 
For this category, the covariant functors $q_*p^*$ of constructible functions 
and  ``twisted" $q_*p^*$ of homologies are  defined, 
and the MacPherson-type natural transformation between these two functors
is constructed. 
In the next section we will construct $C_*^G$ 
throughout applying this construction  to  $\tilde{\phi}_{U,U'}$ and
$\varphi_{U,U'}$. 
\end{rem}

%%%%%%%%%%%%%%%%%%%%%%%%%%%%%%%%%%%%%%%%%%%

\section{Equivariant natural transformation} 
In this section we prove Theorem \ref{C} 
in both contexts of the complex case 
and the quasi-projective case of characteristic $0$ 
(see subsection \ref{mixed}).  
In the latter case, homologies should be read off as Chow homologies $A_*$.  
\subsection{Construction of $C^G_*$}  
For each $U=V-S \in I(G)$ which is non-empty,  the inclusion $U\subset V$ is denoted by $j_U$
and it induces 
$j_U^*:\F^G_{inv}(X\times V) \to \F^G_{inv}(X\times U)$.  
Since $G$ acts freely on $X \times U$  
(hence $X \times U \to X\times_G U$ is a principal bundle), 
any $G$-invariant reduced subscheme $W$ of $X$ has a principal quotient $W
\to W/G$. 
Thus, to $\jeden_W \in \F^G_{inv}(X\times U)$ 
we assign $\jeden_{W/G} \in  \F(X\times_G U)$, 
that actually makes an isomorphism of groups: 
so we identify $\F^G_{inv}(X\times U) =\F(X\times_G U)$.   

As noted in the subsection \ref{mixed},  
we can  apply the (ordinary) MacPherson transformation to the mixed quotient  $X \times_G U$: 
$$C_*: \F(X \times_G U) \to H_*(X \times_G U).$$ 
We denote by  $TU_G$  for short, the vector bundle  
$$X \times_G TU (= X \times_G (U\oplus V)) \to X \times_G U$$ 
and its Chern class  by  $c(TU_G) \in H^*(X \times_G U)$. That is, 
$c(TU_G):=r_U c^G(V)$, 
where $r_U: H^*_G(X) \to H^*(X \times_G U)$ is the canonical projection and 
$c^G(V)$ is the Chern class of the representation $V$. 
Combining the above maps,  we define 
$$T_{U, *}=c(TU_G)^{-1}\frown C_* \circ j_U^*: 
\F^G_{inv}(X \times V)  \to H_*(X \times_G U).$$
Its projection to the truncated homology is also denoted by the same letter.

\begin{lem}\label{com12}
For $U<_* U'$, the following diagram commutes: 
$$
\begin{array}{ccc}
\F^G_{inv}(X \times V) & \stackrel{T_{U, *}}{\longrightarrow} & 
H_{trunc}(X\times_G U)\\
\phi_{U,U'}\downarrow \qquad \quad  & & \qquad \quad \downarrow \varphi_{U,U'}\\
\F^G_{inv}(X \times V')  & 
\stackrel{T_{U', *}}{\longrightarrow} &  
H_{trunc}(X\times_G U')
\end{array}
$$
\end{lem}

\noindent 
\proof:  
We write $T_{U, *}=C_{U, *}\circ j_U^*$  ($C_{U, *}:=c(TU_G)^{-1}\frown
C_*$) for short.
For $U<_*U'$ ($V'=V\oplus V_1$, $U=U-S$, $U'=V'-S'$), 
we take the following diagram, in which 
the left vertical map $p^*$ is $\phi_{U,U'}$, 
the middle vertical map $\iota_*\circ p^*$ is $\tilde{\phi}_{U,U'}$ and 
the right vertical map $\iota_*\circ p^*$ is $\varphi_{U,U'}$: 
$$
\begin{array}{lclcl}
\F^G_{inv}(X \times V) & \overset{j_U^*}{\longrightarrow} &\F(X \times_G U)   
& \overset{C_{U, *}}{\longrightarrow} &  H_{trunc}(X\times_G U)\\
\qquad   \downarrow \; p^* &(1) & \qquad \downarrow\;    p^*  &(2)  & 
\qquad \downarrow \; p^*\\
\F^G_{inv}(X \times V') & 
\overset{j_{U\oplus V_1}^*}{\longrightarrow} &\F(X \times_G (U\oplus V_1))   &
\overset{C_{U\oplus V_1, *}}{\longrightarrow} &  
H_{trunc}(X\times_G (U\oplus V_1))\\
\qquad   \parallel &(3)  &\qquad \downarrow \; \iota_*    &(4)& \qquad
\downarrow \; \iota_*\\
\F^G_{inv}(X \times V') & \underset{j_{U'}^*}{\longrightarrow} &
\F(X\times_G U')   &
\underset{C_{U', *}}{\longrightarrow} &  
H_{trunc}(X\times_G U')
\end{array}
$$
We show that maps of the big square surrounding the diagram commutes. 
First, 
since $\tilde{\phi}_{U,U'}$ is a restriction of $\phi_{U,U'}$ 
as noted in  Remark \ref{Remarkconst}, 
the left half of the big square  ``(1) +(3)" 
 (forgetting the middle arrow $j^*_{U\oplus V_1}$) commutes, 
although (1) commutes but (3) does not. 
Next look at (2). 
For the projection $p=p_{U,U\oplus V_1}: X\times_G (U\oplus V_1) \to
X\times_G U$, 
we let $c(p)$ denote the Chern class of relative tangents of $p$, i.e., 
$c(p)=r_{U\oplus V_1} c^G(V_1) \in H^*(X\times_G (U\oplus V_1))$. 
It then follows from the Verdier-Riemann-Roch formula (\cite{FM},
\cite{Y1}) that 
the following diagram commutes: 
$$
\begin{array}{ccc}
\F(X \times_G U) & \stackrel{C_*}{\longrightarrow} &  H_*(X\times_G U)\\
p^*\downarrow \qquad   & & \qquad \quad \downarrow c(p)\frown p^*\\
\F(X \times_G (U\oplus V_1))   & 
\stackrel{C_*}{\longrightarrow} &  
H_*(X\times_G (U\oplus V_1))
\end{array}
$$
Then for $\alpha \in \F(X \times_G U)$ we have 
\begin{eqnarray*}
p^*C_{U,*}(\alpha) &=& 
p^*\left(c(TU_G)^{-1}\frown C_*(\alpha) \right)   \\
&=&p^*c(TU_G)^{-1}\frown p^*C_*(\alpha)   \\
&=&p^*c(TU_G)^{-1}r_{U\oplus V_1} c^G(V_1)^{-1} 
\left(r_{U\oplus V_1} c^G(V_1)\frown p^*C_*(\alpha) \right)  \\
&=&c(T(U\oplus V_1)_G)^{-1} \frown (c(p) \frown p^*C_*(\alpha) )   \\
&=&c(T(U\oplus V_1)_G)^{-1} \frown C_*(p^*(\alpha)) \\
&=& C_{U\oplus V_1,*}(p^*\alpha). 
\end{eqnarray*}
So the square (2) commutes. 

The remaining part is (4).  
The inclusion $\iota: U\oplus V_1 \to U'$ is not proper 
so we don't use the functoriality of $C_*$, 
but we recall that its complement $U'-U\oplus V_1$ has a sufficiently large 
(complex)  codimension $\ge s$ ($s=\codim  S$ of $U$). 
Therefore, for any $\jeden_W \in \F(X \times_G (U\oplus V_1))$, 
the difference between $c=\iota_*C_*(\jeden_W)$ and $c'=C_*(\iota_*\jeden_W)$ 
has support with $\codim \ge  s$. 
Since $c(T(U\oplus V_1)_G)=\iota^*c(TU'_G)$ is obvious, we have 
$$\iota_*C_{U\oplus V_1, *}(\jeden_W)=c(TU'_G)^{-1}\frown c
\equiv c(TU'_G)^{-1}\frown c' = C_{U', *}(\iota_*\jeden_W)$$
up to the truncated part (of the dimension $>2(n+l'-g)-2s$). 
Thus (4) commutes.  This completes  the proof.  %\qed

\

\begin{definition} \label{CG*}
We define the limit homomorphism 
$$C_*^G:=\lim_{\longrightarrow} T_{U,*}: \F^G(X) \to H_*^G(X), 
\qquad \phi_U(\alpha_U) \mapsto \varphi_U \circ T_{U,*}(\alpha_U). $$
where $\phi_U: \F^G_{inv}(X\times V) \to \F^G(X) $ and 
$\varphi_U: H_{trunc}(X\times_G U) \to H_*^G(X)$ are identification maps. 
\end{definition}

\

\subsection{Proofs of Theorem \ref{C}}
The rest is to show the following properties: 
\begin{enumerate}
\item[(a)] 
$C_*^G:\F^G(X) \to H_*^G(X)$ satisfies the expected naturality: 
\begin{enumerate}
\item[(i)]   
$C_*^G(\alpha+\beta)=C_*^G(\alpha)+C_*^G(\beta)$, $\alpha, \beta \in \F^G(X)$; 
\item[(ii)]  
$C_*^G\circ f_*^G=f_*^G\circ C_*^G$ for proper $G$-morphisms $f: X \to Y$; 
\item[(iii)]  
If $X$ is non-singular, then $C_*^G(\jeden_X^G)=c_*^G(TX)\frown [X]_G$. 
\end{enumerate}
\item[(b)] 
Suppose that for each $U$ we are given a homomorphism 
${DT}_{U,*}: \F^G_{inv}(X\times V) \to H_{trunc}(X\times_G U)$ 
commuting with the structure homomorphisms ($\phi_{U,U'}$ and $\varphi_{U,U'}$) 
such that its inductive limit $D^G_*:\F^G \to H^G_*$ 
satisfying the above properties (i), (ii) and (iii). 
Then $D^G_*$ coincides with our $C^G_*$. 
\end{enumerate}

\proof: 
(a): 
 (i) is trivial.   
 (ii)  follows from the fact that 
$C_*\circ (f\times_G id)_*=(f\times_G id)_*\circ C_*$ 
for the induced proper map $f\times_G id: X\times_G U \to Y\times_G U$. 
 To show (iii)  we recall the normalization property of $C_*$. 
If $X$ is nonsingular, then $X\times_G U$ is also, hence 
$C_*(\jeden_{X\times_G U})=c(T(X\times_G U)) \frown [X\times_G U]$. 
Since $c(T(X\times_G U))=c(TX_G)c(TU_G)$ 
where $c(TX_G)=r_U c^G(TX)$, we have 
$$c(TU_G)^{-1}\frown C_*(\jeden_{X\times_G U})=r_U c^G(TX) \frown
[X\times_G U].$$
Thus, by definitions, 
\begin{eqnarray*}
C_*^G(\jeden_X)&=&C_*^G(\phi_U \jeden_{X\times V}) 
=\varphi_U \circ T_{U,*}(\jeden_{X\times V})  \\
&=& \varphi_U \left( c(TU_G)^{-1}\frown C_*(\jeden_{X\times_G U})  \right)\\
&=&\varphi_U \left( r_U c^G(TX) \frown [X\times_G U] \right)\\
&=& c^G(TX) \frown [X]_G.  
\end{eqnarray*}

\noindent
(b): This is straightforward from the uniqueness of ordinary $C_*$: 
By (iii)  it turns out that ${DT}_{U,*}$ also have 
the natural functoriality and the normalization condition 
up to the truncated homologies. 
Hence by the uniqueness of $C_*$, 
${DT}_{U,*}$ must coincide with our $T_{U,*}$ (up to the truncated homologies). 
Since we can take $U(=V-S)$  so that $\codim S$ is any large number, 
thus  $D^G_*=C^G_*$. 

This completes the proof of Theorem 1.1. %\qed

\subsection{Remark on limit systems}\label{RemarkLim}
It is possible to take some different inductive systems for the definition
of equivariant (co)homology. 
We may replace the structure homomorphisms 
$r_{U',U}:H^*(X\times_G U') \to H^*(X\times_G U)$ 
and $\varphi_{U,U'}:H_{trunc}(X\times_G U) \to H_{trunc}(X\times_G U)$  
 for $U<_* U'$ ($V'=V\oplus V_1$)  by 
$$\tilde{r}_{U',U}:= r_{U'}c^G(V_1)\cdot  r_{U',U}, 
\qquad \tilde{\varphi}_{U,U'}:= r_{U'}c^G(V_1)\frown \varphi_{U,U'}.$$ 
Let $h_G^*(X)$ and $h_*^G(X)$ denote the  limits of these ``twisted"
systems, in a moment. 
This definition of $h_*^G(X)$  is  very  similar to  
the construction of a certain homology theory of a provariety 
(Yokura \cite{Y2}). 

We define $c(T_{BG}) \in h_G^*(X)$ 
to be the projective limit $\{c(TU_G)\}_U$,  which is well-defined: 
$$\tilde{r}_{U',U} c(TU_G)=r_{U'} c^G(V_1)\cdot r_{U',U} c(TU_G)
=r_{U'} c^G(V_1)\cdot  r_{U'}c^G(V)=c(TU'_G).$$
In fact  $c(T_{BG})$ corresponds to  the identity element $1 \in H_G^*(X)$. 
Further we have a group isomorphism 
$c^{-1}:  h_*^G(X) \to H_*^G(X)$ defined by an automorphism of
$H_{trunc}(X\times_G U)$, 
$\xi \mapsto c(TU_G)^{-1}\frown \xi$. 
Then, the above construction of $C^G_*$ is factored as follows: 
$$\F^G(X) \stackrel{\lim\, C_*}{\longrightarrow} h_*^G(X) 
\stackrel{c^{-1}}{\longrightarrow} H_*^G(X).$$
The first map is just the inductive limit of (ordinary) MacPherson
transforms $(C_*)_U$, 
which is well-defined by the VRR 
(this is actually discussed in \cite{Y2}  in a more general setting). 
The second map corrects the twisting of the limit systems, 
that roughly means reducing 
the Chern class factor of 
``horizontal tangents $\pi^*T_{BG}$" of the fibration $\pi: X\times_G EG
\to BG$.

 \subsection{MacPherson's transformation for quotient stacks}\label{Chernq} 
This subsection is a bit isolated from others. 
We work here over $k$ of characteristic $0$ with use of  $A_*$ and assume again 
$G$-schemes  $X$  to be quai-projective  with linearlized $G$-action. 
``The quotient of $X$ via $G$" exists in the category of algebraic stacks, i.e., 
the quotient stack ${\cal X}=[X/G]$ 
associated to the groupoid $G \times X \to X \times X$  (\cite{Vistoli} Appendix, \cite{EG}). 
Then ${\cal X}$  is a category itself, whose objects are  
principal $G$-bundles $E \to T$ 
together with $G$-equivariant morphism $E \to X$ 
and its arrows are morphisms  between principal bundles 
which make the equivariant morphisms to $X$ commute.

In \cite{EG} (Proposition 16) ,  
{\it the integral Chow groups of a quotient stack} ${\cal X}=[X/G]$ 
is introduced as 
$\bar{A}_*({\cal X}):=A^G_{*-g}(X)$  where $g=\dim G$, 
which is independent from the choice of the presentation 
(to avoid  any  confusion,  we put a bar over the letter $A_*$).   
The following lemma is shown in the same manner: 

\begin{lem} 
If $[X/G] \simeq [Y/H]$ as quotient stacks, then $\F^G_{inv}(X) \simeq
\F^H_{inv}(Y)$.  
\end{lem} 

\proof 
Let ${\cal X}:= [X/G]$.  
Since the diagonal of a quotient stack is representative, 
the fibre product $X \times_{\cal X} Y$ is a scheme. 
It has an obvious action of $G \times H$. 
The diagram 
$X \stackrel{p_1}{\leftarrow} X \times_{\cal X} Y
\stackrel{p_2}{\rightarrow} Y$ 
of $G$ and $H$-equivariant projections 
is regarded as the object corresponding to the morphism $Y \to {\cal X}$ 
($p_2$ being the principal $G$-bundle), 
as well the object corresponding to $X \to {\cal X}(\simeq [Y/H])$ 
($p_1$ being the principal $H$-bundle). 
Then via the pullbacks of constructible functions, we have isomorphisms 
$\F^G_{inv}(X) \stackrel{p_1^*}{\rightarrow} 
\F^{G\times H}_{inv}(X \times_{\cal X} Y) \stackrel{p_2^*}{\leftarrow}
\F^H_{inv}(Y)$. 
%\qed

Thus we define the Abelian group of constructible functions for ${\cal
X}=[X/G]$ 
to be  $\bar{\F}_{inv}({\cal X}):=\F^G_{inv}(X)$. 

\begin{rem} \label{remAq}
If $G$ acts on $X$ trivially, $\bar{A}_*({\cal X})$ is identified with
$A_*(X)$. 
As in \cite{EG}, if the action of $G$ is (locally) proper, then 
$\cal X$ becomes an algebraic space 
so the Chow group $A_*({\cal X})$  makes sense, and then  Theorem 3 in
\cite{EG}  says that 
$ A_*({\cal X})\otimes \Q \simeq A^G_{*-g}(X)\otimes \Q \, 
(=\bar{A}_*({\cal X})\otimes \Q)$. 
\end{rem}

Let us consider the following category:  
the objects are quotient stacks  $[X/G]$ of 
quasi-projective varieties $X$ with linearlized actions of some algebraic
groups $G$, 
the arrows are morphisms $\bar{f}: [X/G] \to [Y/H]$, $f: X \to Y$ being
quasi-projective. 

It is easily verified that for this category of quotients,   
$\bar{\F}$ and $\bar{A}_*$ are  covariant functors. 
It follows from Theorem \ref{C} that

\begin{thm}\label{Cq}
For the  category of quotients ${\cal X}$ having presentations $[X/G]$ as
above, 
there is a natural transformation $C_*:\bar{\F}_{inv}({\cal X}) \to
\bar{A}_*({\cal X})$ so that 
for any nonsingular varieties ${\cal X}=X$ (with trivial actions), 
it holds that $C_*(\jeden_X)=c(TX)\frown [X]$. 
\end{thm}

It is also possible to define  $\F([X/G])$  as the inductive limit of 
$\F_{inv}([X\times V/G])$ and extend $C_*$ to the group 
(cf.   Corollary \ref{VRRcor} below).  

As noted in  \ref{mixed},  
presumably the above theorem would be stated 
for general quotient stacks (i.e., without quasi-projectiveness)  
by generalizing Kennedy's formulation to algebraic spaces. 
Not only for quotient stacks, 
also for (general)  algebraic stacks, MacPherson's transformation $C_*$ is 
expected, 
but the author does not know how to manage it.

%%%%%%%%%%%%%%%%%%%%%%%%%%%%%%%%%%%%%%%%%%%
\section{Some properties of $C^G_*$}
\subsection{Components of $C^G_*(\jeden_X)$}\label{degree} 
It is clear by definition that 
$C^G_*(\jeden_X)$ consists of $C^G_i(\jeden_X) \in H^G_{2i}(X)$, 
which is called  {\it the $i$-th component}, 
furthermore  
$C^G_i(\jeden_X)$ is always trivial for any negative $i$ 
 (e.g., see Corollary \ref{negativecomp}).  

We remark about the lowest $0$-th component and the top $n$-th one of
$C^G_*(\jeden_X)$ 
for a projective $G$-variety $X$ of (equi)dimension $n$. 

Let $f: X\to \{pt\}$, the pointed map. 
For any invariant constructible function $\alpha \in \F^G_{inv}(X)$ 
{\it the degree} of the $0$-th component $G_0(\alpha)$ is defined to be the
number 
$f^G_*C^G_0(\alpha) \in \F^G_{inv}(pt)=\Z$. It is easily seen that 
$$\mbox{the degree of $C^G_0(\alpha)$} = \int_X \alpha \;\; \in \; \Z.$$
In fact, for $\alpha = \jeden_X$, 
$$f^G_*C^G_0(\jeden_X)\stackrel{(1)}{=}C^G_0 f^G_*(\jeden_X)
\stackrel{(2)}{=}\chi(X) C^G_0(\jeden_{pt})\stackrel{(3)}{=}\chi(X) \;
[pt]_G. $$
Here (1) is by the naturality, (2) comes from the linearity and the fact that 
$f^G_*(\jeden_X)=\int_X \jeden_X = \chi(X)$ as seen before, 
and (3) is the normalization condition 
$C^G_*(\jeden_{pt}) =c^G(Tpt)\frown [pt]_G=[pt]_G$ 
(corresponding to $1 \in  H^*(BG)$).

As to the top component, recall the ordinary case: $C_n(\jeden_X)=[X] \in
H_{2n}(X)$.  
In our equivariant setting, 
$C^G_*(\jeden_X)$  is the limit of 
$c(TU_G)^{-1}\frown C_*(\jeden_{X\times_G U})$ 
whose top component is the fundamental class $[X\times_G U]$, thus 
$C^G_n(\jeden_X) =[X]_G$.

\subsection{Change of groups}\label{changegroup}
Let $G_0$ be the closed subgroup of $G$, and $X$ a $G_0$-variety of
dimension $n$. 
Then $X\times_{G_0} G$ becomes a $G$-variety of dimension $n+k$ ($k=\dim
G/G_0$).
The $G$-constructible functions and equivariant homology of the mixed space 
are identified as follows: 

\begin{prop}
There are canonical isomorphisms so that the following diagram commutes: 
$$
\begin{array}{ccccc}
\F^{G}(X\times_{G_0} G) 
& \stackrel{C^{G}_*}{\longrightarrow} & H^{G}_{2(n+k)-*}(X\times_{G_0} G) 
& \stackrel{\cap \mbox{\tiny fund}}{\longleftarrow} & H_G^*(X\times_{G_0} G) \\
 \simeq\downarrow \;\; & & \;\; \downarrow \simeq 
& & \;\; \downarrow \simeq \\
\F^{G_0}(X) & \stackrel{C^{G_0}_*}{\longrightarrow} & H^{G_0}_{2n-*}(X) 
& \stackrel{\cap \mbox{\tiny fund}}{\longleftarrow}  & H_{G_0}^*(X) 
\end{array}
$$
\end{prop}

\proof  
Take $U \in I(G)$, an open set of a representation $V$  of $G$ over which
the action is free. 
Of course, $G_0$ acts freely on $U$, i.e., $U \in I(G_0)$.   
We denote simply a point of $(X \times_{G_0} G) \times_G U$ by $[[x,a],u]$ 
so that $[[x,a],u]=[[h.x, ha], u]=[[h.x,hag^{-1}],g.u]$ for $h \in G_0$ 
and $g \in G$. 
Then,  to each $[[x,1], u]$, we assign $[x,u]$ to get 
a (well-defined)  isomorphism 
$(X\times_{G_0} G) \times_G U \stackrel{\sim}{\to}  X \times_{G_0} U$.  
Also we can see an one-to-one correspondence between invariant subvarieties 
of  $(X\times_{G_0} G) \times V$ and the one of $X\times V$. 
By these identifications and the construction of $C^G_*$, the claim
follows. %\qed

\subsection{Cross products} 
In (ordinary) MacPherson's theory, the cross product formula is known (\cite{Kw2}): 
for $\alpha \in \F(X)$ and $\beta \in \F(Y)$, 
$\alpha \times \beta\in \F(X\times Y)$  is defined to be 
$\alpha \times \beta(x,y) := \alpha(x) \cdot \beta(y)$, and it then holds that 
$C_*(\alpha \times \beta) = C_*(\alpha) \times C_*(\beta)$ 
where $\times$ in the right hand side means the homology cross product. 

For $G$-varieties $X$ and $Y$, the cross product 
$$H^G_i(X) \otimes H^G_j(Y) \to H^G_{i+j}(X\times Y), \quad 
(\xi, \xi') \mapsto \xi \times \xi'$$ 
is well-defined (Def-Prop. 2 of \cite{EG}). In fact, 
for any $U=V-S$ and $U'=V'-S'$ in $I(G)$, 
we set the isomorphism 
$s:(X\times U)\times (Y\times U') \to  X\times Y \times (U\oplus U')$ 
by $(x,u,y,u') \mapsto (x,y,u\oplus u')$  
(also the isomorphism between the quotients is denoted by $s$), 
and then the limit of 
$$H_{trunc}(X\times_G U) \oplus H_{trunc}(Y\times_G U')  
\to H_{trunc}((X\times Y) \times_G (U\oplus U')), $$
 $(c, c') \mapsto s_*(c \times c') $, 
gives the above cross  product. 
In a similar way,  we define the exterior product of $G$-constructible
functions, 
$$\F^G_{inv}(X\times U) \otimes \F^G_{inv}(Y\times U') 
\to \F^G_{inv}((X\times Y) \times_G (U\oplus U')), $$
$(\alpha, \beta) \mapsto s_*(\alpha \times \beta)$. 
By the (ordinary)  theory,  we see that 
$C_*s_*(\alpha \times \beta)=s_*C_*(\alpha \times \beta)
=s_*(C_*(\alpha) \times C_*(\beta))$. 
Multiplying both sides by inverse Chern class factors ($c(TU_G)^{-1}$ etc) 
and then taking the limits, 
we have 
$$C^G_*(\alpha \times \beta) = C^G_*(\alpha) \times C^G_*(\beta) 
\quad (\alpha \in \F^G(X), \; \beta \in \F^G(Y)).$$

\subsection{Equivariant Verdier-Riemann-Roch} 
We can formulate 
an equivariant version of the VRR formula for smooth morphisms (\cite{Y1}).  

\begin{thm}\label{VRR} 
Let $f: X \to Y$ be a  $G$-equivariant smooth morphism 
with a $G$-equivariant relative tangent bundle $\nu$ of dimension $m$. 
Then  the following diagram commutes: 
$$
\begin{array}{ccc}
\F^G(Y)  & \stackrel{C^G_*}{\longrightarrow} & H^G_*(Y)\\
f^*\;\downarrow \;\;& & \;\; \downarrow \; f^{**} \\
\F^G(X)  & \stackrel{C^G_*}{\longrightarrow} & H^G_{*+2m}(X)
\end{array}
$$
where $f^{**}$ is defined by $f^{**}(c):=c^G(\nu) \frown f^*(c)$. 
\end{thm}

\proof 
Given a $f:X \to Y$ as above, then 
for any $U$,   we set $f':=f\times_G id :X\times_G U \to Y\times_G U$, 
which is an smooth morphism with 
the relative bundle $\nu_G:=\nu\times_G U$ over $X\times_G U$. 
We apply the ordinary VRR formula to this $f'$, 
and then we have 
$$
\begin{array}{ccc}
\F(Y\times_G U)  & \stackrel{C_{U, *}}{\longrightarrow} & H_*(Y\times_G U)\\
(f')^*\;\downarrow \;\; \; & & \;\;\; \downarrow \; (f')^{**} \\
\F(X\times_G U)  & \stackrel{C_{U, *}}{\longrightarrow} & H_{*+2m}(X\times_G U)
\end{array}
$$
where $C_{U, *}=c(TU_G)^{-1}\frown C_*$ and $(f')^{**}= c(\nu_G)\frown (f')^*$. 
This commutative diagram for $U$ and the one for $U'$ with $U<_*U'$ 
are compatible with $\tilde{\phi}_{U,U'}$ and $\varphi_{U,U'}$, 
so we see the claim  
similarly  as in the proof of Lemma \ref{com12}. 
%\qed

\begin{cor}\label{VRRcor} 
The limit homomorphism 
$$\underset{\longrightarrow}{\lim}\; (C^G_*\circ \phi_0): 
\F^G(X)=\underset{\longrightarrow}{\lim}\; \F^G_{inv}(X\times V)  \to 
\underset{\longrightarrow}{\lim}\;H^G_*(X\times V) $$
coincides with the composition  $p^{**}\circ C^G_*$ of 
the transformation $C^G_*:\F^G(X) \to H^G_*(X)$ 
and the twisted identification 
$p^{**}: H^G_*(X) \stackrel{\sim}{\to} H^G_{*+2\dim V}(X\times V)$, 
$p^{**}=c^G(V)\frown p^*$, where $p$ is the projection $X \times V \to X$. 
\end{cor}

\proof 
The claim is shown by the following diagram whose right square is a simple VRR 
$$
\begin{array}{ccccc}
\F^G_{inv}((X\times V) \times V) & \overset{\phi_{V}}{\longrightarrow} & 
\F^G(X\times V)  & \overset{C^G_*}{\longrightarrow} & 
H^G_{*+2l}(X\times V) \\
\quad\uparrow (p\times id)^*& & \quad \uparrow\; p^* & & \simeq \; 
\uparrow \; p^{**}\\
\F^G_{inv}(X\times V) & \underset{\phi_V}{\longrightarrow} 
&\F^G(X)  & \underset{C^G_*}{\longrightarrow} & H_*^G(X)
\end{array}
$$
where first two $\phi_V$'s are identification maps of the inductive limits. 
Clearly,  $\phi_{V}\circ (p\times id)^*$ equals the canonical inclusion 
$\phi_{0}:\F^G_{inv}(X\times V) \subset \F^G(X\times V)$. 
%\qed

%%%%%%%%%%%%%%%%%%%%%%%%%%%%%%%%%%%%%%%%%%%
\section{Equivariant Chern-Mather class} 
\subsection{Construction of $C^M_G$} 
First of all, we recall that 
there are  two key factors in  (ordinary)  MacPherson's Chern class,  
{\it the Chern-Mather class} and 
{\it the local Euler obstruction}, 
see \cite{Mac} for the detail. 
The local Euler obstruction $Eu_X $ of a variety $X$ is 
a constructible function of $X$ 
assigning to $x \in X$  an local invariant of the germ $(X,x)$, 
which is defined using obstruction theory (through the transcendental
method)  in \cite{Mac}, 
see also \cite{Gon} for the purely algebraic definition. 
It defines an isomorphism
$Eu: Z(X) \to \F(X)$ sending algebraic cycles $\sum n_k W_k$ 
to $\sum n_k Eu_{W_k}$, 
and $C_*$ is constructed so that 
$C_*(\sum n_k Eu_{W_k})=\sum n_k (i_k)_*C^{M}(W_k)$ 
where $C^M$ means the Chern-Mather class, 
$i_k:W_k \to X$ is the inclusion map and $(i_k)_*$ is the induced homomorphism 
of homology groups. 

Almost automatically, 
this relation among $C_*$, $C^M$ and $Eu$  is lifted into the equivariant
version. 
Let $X$ be a $G$-reduced scheme of equidimension $n$ and 
assumed to be $G$-embeddable into some $G$-nonsingular variety, say $M$. 
As the same in the ordinary case without actions, 
we have {\it the $G$-equivariant Nash blow-up} of $\nu:\widehat{X} \to X$: Here 
 $\widehat X$  is given as  the closure of the regular part $X_{Reg}$ in 
the Grassmanian $Gr_n(TM)$ of $n$-planes of $TM$ on which $G$ acts naturally, 
and $\nu:\widehat X \to X$  is the natural projection, 
which is an equivariant proper morphism. 
$$
\begin{array}{ccc}
\widehat X & \subset & Gr_n(TM)\\
\nu\;\downarrow\quad  & & \downarrow\\
X & \subset & M
\end{array}
$$
Let $\widehat {TX}$ be the $G$-Nash tangent bundle over $\widehat{X}$, 
that is the restriction over  $\widehat{X}$ 
of the tautological $G$-vector bundle of the Grassmanian. 
Then, 
we define {\it the $G$-equivariant  Chern-Mather class of $X$}  to be 
$$ C^M_G(X) := \nu_*^G(c^G(\widehat {TX})\frown [\widehat X]_G) \in H^G_*(X).$$
Note that $\nu$  is actually made 
by the local embedding of $X$ and the gluing  process, 
so this construction is independent from the choice of the ambient space 
(cf. \cite{Mac}, \cite{Ken}). 

For a $G$-variety $X$, 
let $Z^G_{inv}(X)$ denote the subgroup of $Z(X)$ consisting of 
$G$-invariant algebraic cycles, 
and let $C^M_G: Z^G_{inv}(X) \to H^G_*(X)$ denote the map sending invariant
cycles 
$\sum n_k W_k$ to $\sum n_k (i_k)^G_*C^M_G(W_k)$ 
where each $W_k$ is a $G$-invariant reduced closed subscheme of $X$. 
We write the canonical inclusion by $\phi_0: \F^G_{inv}(X) \subset \F^G(X)$. 

\begin{prop}\label{CM}  Let $X$ be a $G$-variety. Then,  \\
(1)  $Eu_X$ is $G$-invariant and $Eu: Z^G_{inv}(X) \to \F^G_{inv}(X)$ is an
isomorphism; \\
(2)  The following diagram commutes: 
$$
\begin{array}{rcl}
Eu\;\; &  Z^G_{inv}(X) & \;\; C^M_G\\
\overset{\simeq}{\swarrow}& & \searrow  \\
\F^G_{inv}(X)& \longrightarrow & H^G_{*}(X)\\
& C^G_*\circ \phi_0 &
\end{array}
$$
\end{prop}

\proof 
(1) For any  reduced $G$-invariant subscheme  $W$ and for any $g\in G$,  
$g: (W,x) \to (W, g(x))$ is isomorphic, so 
$Eu_W(x)=Eu_W(g(x))$, i.e., $Eu_W \in \F^G_{inv}(X)$.  
Since $Eu:Z(X) \to \F(X)$ is an isomorphism, its restriction to
$Z^G_{inv}(X)$ is 
an isomorphism into $\F^G_{inv}(X)$.  
The surjectivity is shown similarly as in \cite{Mac}  by the induction of
dimension using 
the fact that the singular loci of a $G$-invariant variety is also
$G$-invariant. 

\noindent 
(2) 
Let $W$ be an invariant reduced subscheme of $X$. 
First, we see 
$$Eu_W(x)=Eu_{W\times U}(x,u)=Eu_{W \times_G U}([x,u]).$$
This follows from three facts: 
$Eu_{W\times W'}(x,y)=Eu_W(x)Eu_{W'}(y)$, 
$Eu_W(x)=1$ if $(W,x)$ is nonsingular, and 
$W\times U \to W\times_G U$ is a principal $G$-bundle. 
Then, the map 
$$ j_U^*\circ \phi_{0,U}:\F^G_{inv}(X) \to \F^G_{inv}(X\times V)  \to
\F(X\times_G U)$$
sends $Eu_W$  to $Eu(W\times_G U)$. 

Second,  let $\nu: \widehat{W} \to W$ be 
the $G$-Nash blow-up of $W$, and set $W_G:=W\times_G U$. 
Since $W_G \to U/G$ is a bundle with fibre $W$ over the smooth base space,  
it turns out that $\nu'=\nu \times_G id: \widehat{W}\times_G U \to W_G$ 
gives the (ordinary) Nash blow-up of $W_G$, i.e., 
$\widehat{W_G}=\widehat{W}\times_G U$.
Then the (ordinary) Nash tangent bundle $\widehat{TW_G}$ 
splits to two factors coming from 
the $G$-Nash tangent bundle $\widehat{TW}$ and $TU$. 
By these facts and the fact that $C^M=C_*\circ Eu$ as mentioned, 
we have 
\begin{eqnarray*}
C_*\circ j_U^*\circ \phi_{0,U}(Eu_W) 
& =& C_*\circ Eu(W\times_G U) =C_*\circ Eu(W_G)  \\
&=& C^M(W_G) = (\nu')_*(c(\widehat{TW_G})\frown [\widehat{W_G}])\\
&=&(\nu')_*c((\widehat{TW}\times_G U)\oplus(\widehat{W}\times_G TU))
\frown  [\widehat{W_G}])\\
&=&c(TU_G)\frown ((\nu')_*(c(\widehat{TW}\times_G U)\frown [\widehat{W}_G]) ).
\end{eqnarray*}
Multiply both sides by $c(TU_G)^{-1}$ and take the limit, then 
we have $C^G_*\circ \phi_0 (Eu_W) = C^M_G(W)$.  %\qed

\

By its definition the Mather class $C^M_G(X)$  
has no negative-dimensional component of $H^G_*(X)$. Thus it follows that 

\begin{cor} \label{negativecomp} 
For any $G$-invariant constructible function $\alpha \in \F^G_{inv}(X)$, 
each negative-dimensional component 
$C^G_i(\alpha)$ ($i<0$) is trivial ($\in H^G_{2i}(X)$). 
\end{cor}

This elementary fact is not very clear in the functorial definition of
$C^G_*$. 
If we assume the existence of $G$-equivariant desingularizations
of $X$ being considered  
(or assume that the quotient $X/G$ becomes a variety), 
then the above corollary  immediately follows from 
functorial and normalization properties of $C^G_*$ 

\

Next, in the diagram of Proposition \ref{CM} let us replace $X$   
by $X\times V$ with a representation $V$ and take the limit of groups. Then
we have 

\begin{cor} \label{CMcor}
Throughout the twisted identification via $p^{**}$ as in Corollary
\ref{VRRcor}, 
$C^G_*$ is  factored  by  the  limit map of $Eu^{-1}$ and the one of
$C^M_G$, i.e., 
$p^{**}\circ C^G_*=\underset{\to}{\lim} \, C^M_G \circ \underset{\to}{\lim}
\, Eu^{-1}$: 
$$
\begin{array}{rcl}
& \underset{\longrightarrow}{\lim}\; Z^G_{inv}(X\times V) & \\
 \simeq \; \nearrow &&  \searrow  \\
\underset{\longrightarrow}{\lim}\; \F^G_{inv}(X\times V)& 
\longrightarrow
&\underset{\longrightarrow}{\lim}\;H^G_{*+2l}(X\times V) \\
| |  \qquad \quad & & \qquad  p^{**} \; \uparrow \; \; \simeq \\
\F^G(X) \quad & \underset{C^G_*}{\longrightarrow} & \qquad  H^G_*(X)
\end{array}
$$
\end{cor}

\begin{rem}\label{ChernM} 
$C^M_G(X)$ is a kind of 
 the Mather  class with respect to  the ``fibrewise Nash-blowing up" of 
$X\times_G U \to U/G$ 
as seen in the proof of Proposition \ref{CM}.  
 This explains again the reason that  $c(TU_G)^{-1}$ appears in the
definition of $T_{U,*}$ 
and  the twisted identification appears 
in Corollaries  \ref{VRRcor} and  \ref{CMcor}   (cf.   Remark \ref{RemarkLim}). 
We may have started to define  the $G$-equivariant MacPherson transformation 
$C^G_*:\F^G_{inv} \to H^G_*$  by $C^M_G\circ Eu^{-1}$ as described above 
(or using a certain analogy  to Lagrangian cycles),  
but it would become the same thing as we have done.  
To prove the functoriality of pushforwards would be harder than our simple
argument. 
\end{rem}

\begin{rem} \label{ChernF} 
As another story of Chern classes for singular spaces, there is 
so-called {\it the Fulton canonical class} $C^F(X)$ (Example 4.2.6 \cite{F}): 
Let $s(X,M)$ denote the Segre covariant class of a closed subscheme $X$  
in a non-singular variety $M$, 
and then $C^F(X)$ is defined to be $c(TM|_X) \frown s(X,M)$. 
It turns out that $C^F(X)$ depends only on $X$, not on the choice of the
ambient space $M$. 
In the equivariant setting, 
it is also possible to make the $G$-version of Fulton's class $C^F_G(X) \in
H^G_*(X)$  
for a $G$-scheme $X$ using an equivariant blowing-up 
by the $G$-invariant ideal sheaf defining $X$ in $M$. 
Then  $C^{SM}_G$, $C^M_G$ and $C^F_G$ 
would be ``basic" equivariant Chern homology classes like as the ordinary case, 
and they would have particular interests in singularity theory.  
\end{rem}

\subsection{Restriction: ordinary transformation without $G$-action} 
Take a fibre of $X\times_G U \to U/G$,  
and denote the inclusion of the fibre by  $i:X \subset X\times_G U$ in
abusing the notation. 
Then we have a homomorphism $i^*: H^G_*(X) \to H_*(X)$, that is independent from 
the choice of fibre. 
Through $i^*$,  our $C^G_*$  recovers the ordinary MacPherson
transformation $C_*$. 

\begin{prop}\label{Restriction}
We have the commutative diagram: 
$$
\begin{array}{ccccc}
  \F^G_{inv}(X) &\overset{Eu^{-1}}{\to}  &  Z^G_{inv}(X)  &
\overset{C^M_G}{\to} & H^G_*(X)  \\
\cap && \cap && \;\;\downarrow \, i^* \\
\F(X) & \overset{Eu^{-1}}{\to} &  Z(X)  & \overset{C^M}{\to}  & H_*(X) 
\end{array}
$$
\end{prop}

\proof 
The left square means the restriction of $Eu^{-1}$ to $\F^G_{inv}(X)$, so
it commutes. 
As seen in the proof of Proposition \ref{CM}, 
the restriction of  $\nu \times_G id: \widehat{X} \times_G U \to X\times_G
U$  to a fibre $i(X)$ 
 is isomorphic to the ordinary Nash blowing-up $\nu:\widehat{X}  \to X$,  thus 
by the definitions of $C^M$ and $C^M_G$, the right square commutes. %\qed

%%%%%%%%%%%%%%%%%%%%%%%%%%%%%%%%%%%%%%%%%%%
 \section{Quotient Chern  classes}
\subsection{Canonical constructible functions} 
We introduce naturally defined {\it canonical constructible functions} 
and the corresponding equivariant Chern classes for a $G$-variety. 

We again abuse the (set-theoretic) notation for simplicity, 
and we denote 
$$\Psi:  G\times X \to X \times X, \; \;\Psi(g,x)= (g.x, x), $$
$$Z :=\Psi^{-1}\mit\Delta_X=\{(g,x)|g.x=x\}. $$ 
The projections of $Z$ to factors are denoted by   
$G \stackrel{p_1}{\leftarrow} Z \stackrel{q_1}{\rightarrow} X$. 
Note that $p_1^{-1}(g)=\{x\in X|g.x=x\}=:X^g$ 
(the fixed point set of the map $g: X \to X$ ($x \mapsto g.x$)) and 
$q_1^{-1}(x)=\{g\in G|gx=x\}=:Stab_x(G)$  (the stabilizer subgroup of  $x$). 
Then we can define a $G$-invariant constructible function on $X$ 
$$\alpha_{X/G}^{(1)} := (q_1)_*\jeden_{Z_{red}}  \; \in \; \F^G_{inv}(X).$$ 
The invariance is clear. In particular,  for $x \in X$ 
$$\alpha_{X/G}^{(1)} (x) 
=\int_{G\times \{x\}} \jeden_{Z_{red}}=\chi(Stab_x(G)_{red}).$$

The adjoint action of $G$ on itself and the preimage of the diagonal are
denoted by 
$$\Phi:  G\times G \to G \times G, \;\; \Phi(h, g)=(hgh^{-1},g), $$
$$Com(G)=\Phi^{-1}\mit\Delta_G=\{(h,g)|gh=hg\},$$ 
 and the projection $r_2:Com(G) \to G$ ($r_2(h,g)=g$). 
Of course, the fibre $r_2^{-1}(g)$ is isomorphic to the centralizer subgroup
$C(g)$ of $g$. 
Let 
$$Z^{(2)}:=\{ \; (h, g, x) \in G \times G\times X \; | \; gh=hg, \; g.x=x,
\; h.x=x\; \}$$
 and the natural projections 
$p_2:Z^{(2)} \to Com(G)$ ($(h,g,x) \mapsto (h,g)$) and $q_2:Z^{(2)}\to Z$
($(h,g,x) \mapsto (g,x)$). 
Then for a point $(g,x) \in Z$, the fibre $(q_2)^{-1}(g,x)$ is isomorphic
to $Stab_x(G) \cap C(g)$ 
(=$Stab_x(C(g))$ for the action of $C(g)$ on $X^g$).  
It is also obvious that for $(h,g) \in Com(G)$, $(p_2)^{-1}(h,g)$ 
corresponds to $X^{\{h, g\}}:=X^h \cap X^g$, 
the intersection of fixed point sets.  
We set $\pi_2:= q_1\circ q_2:Z^{(2)}\to X$ and define 
$$\alpha_{X/G}^{(2)}:= (\pi_2)_*\jeden_{Z^{(2)}_{red}}  \in \F^G_{inv}(X).$$
This procedure can be continued  by introducing 
the set of mutually commuting $k$-tuple of elements of $G$, 
$Com(G;k):=\{(g_k, \cdots, g_1)| g_ig_j=g_jg_i\}$ (cf.  \cite{BF}) 
and the correspondence 
$$Z^{(k)}:=\{ \; (g_k, \cdots, g_1,x) \in Com(G;k)\times X \; |\;  g_i \in
Stab_x(G)\;  \}$$
(e.g., $Com(G;1)=G$, $Com(G;2)=Com(G)$, $Z^{(0)}=X$, $Z^{(1)}=Z$). 
The natural projections to  $Com(G;k)$ and $Z^{(k-1)}$ 
are  defined, thus we have the following commutative diagrams  
$$
\begin{array}{ccccccc}
Z^{(k)} & \stackrel{q_k}{\to} &Z^{(k-1)} & \stackrel{q_{k-1}}{\to}  \cdots
\stackrel{q_2}{\to}  
& Z  & \stackrel{q_1}{\to} & X \\
\mapdL{p_k\; \;\;  }\;  & & \mapdL{p_{k-1}\; \;\;  }\;  
&  &\mapdL{p_1\; \; \;  }  &&  \mapdL{p_0\; \;
\;  }  \\
Com(G;k)  &\stackrel{r_k}{\to}  & Com(G;k-1) 
&\stackrel{r_{k-1}}{\to}  \cdots \stackrel{r_3}{\to} 
&  G &\to & \{e\}
\end{array}
$$
The preimage $p_k^{-1}(g_k,\cdots, g_1)$ is 
the set of simultaneous fixed points, denoted by $X^{\{g_k,\cdots, g_1\}}$. 
Let  $\pi_0=id_X$ and $\pi_k:=q_1\circ \cdots \circ q_k: Z^{(k)}  \to  X$. 
Then  
the $k$-th {\it canonical constructible functions} of $X$ is defined by  
$\alpha_{X/G}^{(0)}:=\jeden_X$ and 
$$\alpha_{X/G}^{(k)}:=(\pi_k)_*\jeden_{Z^{(k)}_{red}}\quad \in \;
\F^G_{inv}(X). $$

Applying our natural transformation to these canonical
constructible functions, 
we obtain  a sequence of {\it integral} Chern classes 
$C_*(\alpha_{X/G}^{(k)})$ in $H^G_*(X)$, 
say {\it canonical quotient Chern classes} for $X/G$.  
Furthermore, 
if $\chi(G)\not=0$, then 
 {\it rational canonical constructible functions} are defined by 
$$\jeden_{X/G,orb}^{(k)}:=\frac{\alpha_{X/G}^{(k)}}{\chi(G)} \quad 
\in\quad \F^G_{inv}(X) \otimes \Q$$ 
as well {\it rational} quotient Chern classes  
$C_*(\jeden_{{\cal X}, orb}^{(k)} )$ in $H^G(X) \otimes \Q$. 
For the importance of the case $k=1, 2$, we write as 
$\jeden_{X/G,quo}:=\jeden_{X/G,orb}^{(1)}$ and
$\jeden_{X/G,orb}:=\jeden_{X/G,orb}^{(2)}$. 

Note again that the canonical quotient Chern classes are defined for
arbitrary actions of $G$ 
($\dim G \ge 0$). 
In the next subsection we will see these classes in a typical case of a
finite group action.

\subsection{Orbifold Chern class}  
Let $G$ be a finite group and $X$ a possibly singular quasi-projective
variety with $G$-action.  
We also assume that $X/G$ is a variety and the stabilizer subgroups are
reduced.  
Combining Proposition \ref{Restriction} with the pushforward property of 
the ordinary $C_*$  for the projection $\pi: X \to X/G$, 
we have the commutative diagram: 
$$
\begin{array}{ccccc}
 \F^G_{inv}(X) & \overset{C^G_*\circ\phi_0}{\longrightarrow} & H^G_*(X)  \\
 \pi_*\, \downarrow \;\;\;   & & \; \;  \downarrow \,  \pi_*\circ i^*   \\
 \F(X/G) & \overset{C_*}{\longrightarrow} & H_*(X/G)  
\end{array}
$$
Since $g=\dim G=0$, 
the map $\pi_*\circ i^*$, for nonnegative dimension $*\ge 0$,  
is an isomorphism after tensored by $\Q$, 
which coincides with the isomorphism mentioned in Remark \ref{remAq} 
(i.e.,  Theorem 3 in \cite{EG}).  
Note that the first horizontal arrow is just $C_*$ given 
for the quotient stack $[X/G]$ in Theorem \ref{Cq}, 
and the above diagram says the simple fact that if $X/G$ is a variety, 
the above $C_*$  is identified with 
ordinary MacPherson's transformation (the lower horizontal arrow) within rational coefficients. 

An elementary computation shows that 
$$\jeden_{X/G, quo} =\frac{1}{|G|}\sum_{g}  \jeden_{X^g}, \quad 
\jeden_{X/G, orb} =\frac{1}{|G|}\sum_{gh=hg}   \jeden_{X^{\{g,h\}}}$$
in $\F^G_{inv}(X)\otimes \Q$, where the first sum runs over all $g \in G$ and 
the second runs over all $(h,g) \in Com(G)$: In fact, 
\begin{eqnarray*}
\alpha_{X/G}^{(1)}(x)&=&|Stab_x(G)|=\sum \jeden_{X^g}(x); \\
\alpha_{X/G}^{(2)}(x)&=&\sharp\; (\pi^{(2)})^{-1}(x)=\sharp\;
\{(g,h)|gh=hg, g.x=x, h.x=x\}  \\
&=& \sum   \jeden_{X^{\{g,h\}}}(x). 
\end{eqnarray*}
We call $C^G_*(\jeden_{X/G, quo})$ {\it the Chern class of $X/G$} and 
$C^G_*(\jeden_{X/G, orb})$ {\it the orbifold Chern class of $X/G$}, in
abusing words. 
Note that if $X$ is nonsingular, 
we can take the dual of these classes in  $H_G^*(X)\otimes \Q$.  

\begin{prop} It holds that $\pi_*(\jeden_{X/G, quo})=\jeden_{X/G}$ and 
$$\pi_*i^*\, C^G_*(\jeden_{X/G, quo})=C^{SM}_*(X/G) 
\quad \in \; \; H_*(X/G;\Q) $$
where $i^*$ is the  map induced by a restriction $i$ given in 4.3 
and $\pi_*$ is the pushforward induced by the projection $\pi$.  
\end{prop}

\proof  
This is also an elementary computation. Since 
$|G|=\sharp\, G/Stab_x(G)\cdot |Stab_x(G)|= \sharp\, G.x \cdot |Stab_x(G)|$, 
and $|Stab_{x'}(G)|=|Stab_x(G)|$ (for $x' \in G.x$), 
we have that for any $[x] \in X/G$ 
\begin{eqnarray*}
\pi_*(\jeden_{X/G, quo}) ([x])&=&\int_{G.x} \jeden_{X/G, quo} 
= \sum_{x'\in G.x} \jeden_{X/G, quo}(x') \\
&=& \frac{1}{|G|}\sum_{x'\in G.x} \;\sum_{g\in G} \; \jeden_{X^g}(x')=1.
\end{eqnarray*}
The second equality follows from the commutative diagram as described
above.  %\qed

\begin{cor}
It holds that 
\begin{eqnarray*}
\chi(X/G)&=& \mbox{\rm degree of $0$-th component}\; \; C^G_0(\jeden_{X/G, quo})  \\
\left( \right.&=&\int_X \jeden_{X/G, quo} 
=\frac{1}{|G|}\sum_{g } \; \chi(X^g) \quad \left.\right). 
\end{eqnarray*}
\end{cor}

\proof  
This follows from the above proposition and the subsection \ref{degree}. 
Let us take  $f: X/G \to \{pt\}$ and $f\circ \pi: X \to \{pt\}$. 
Then by using the above proposition, we have 
$$
\int_X \jeden_{X/G, quo} =  (f\circ\pi)_*(\jeden_{X/G, quo})  = 
f_*\pi_*(\jeden_{X/G, quo}) 
=   f_*\jeden_{X/G}   = \chi(X/G).  %\quad \mbox{\qed}
$$

Next, let us look at  {\it the  orbifold Euler characteristics} of the
quotient $X/G$, then 
the same type equalities hold: For a possibly singular $G$-variety $X$, 
\begin{eqnarray*}
\chi(X; G) &=& \mbox{degree of $0$-th component}\; \; C^G_0(\jeden_{X/G,
orb})  \\
\left( \right.&=&\int_X \jeden_{X/G, orb}= \frac{1}{|G|}
\sum_{gh=hg} \; \chi(X^{\{g,h\}}) \quad \left.\right). 
\end{eqnarray*}
We may expect that $\chi(X; G)$ is related to 
the Euler characteristics of some certain desingularizations of $X/G$ (cf. \cite{HH}). 
Here is an optimistic conjecture: 

\begin{conj}  Let $G$ be a finite group, or more generally a linear
algebraic group.  
Let $X \to Y:=X/G$ be the quotient. 
Then there is a $G$-variety $\tilde{X}$ and a proper $G$-morphism
$f:\tilde{X}  \to X$ 
so that  $f^G_*(\jeden_{\tilde{X}/G, quo}) = \jeden_{X/G, orb}$. In particular, 
$f^G_*C^G_*(\jeden_{\tilde{X}/G, quo}) = C^G_*(\jeden_{X/G, orb})$ and 
their  $0$-th degrees $\chi(\tilde{X}/G)=\chi(X;G)$  
\end{conj}

Also we obtain a certain sequence of 
the generalized orbifold Chern classes $C^G_*(\jeden^{(k)}_{X/G, orb})$,  
whose $0$-th component is just 
the generalized orbifold Euler characteristics defined in \cite{BF}: 
$$\mbox{degree of  $C^G_0(\jeden^{(k)}_{X/G, orb})$} 
= \frac{1}{|G|} \sum_{Com(G;k)} \chi(X^{\{g_1, \cdots, g_k\}}),$$
where $X^{\{g_1, \cdots, g_k\}}$ is the simultaneously fixed point set. 
As a particularly interesting case, 
the Chern class of symmetric products 
 ($S^nX=X^n/S_n$, $S_n$ being the $n$-th symmetry group) 
is studied in \cite{O2}, in which 
{\it the generating functions of orbifold  Chern classes} of $S^nX$ is
obtained. 

Our principle is that  
{\it certain kinds of formulas on Euler characteristics 
should admit the Chern class version}. 
From this viewpoint, 
``the constructible function-description" is seemingly very useful 
and very natural. Also we have another advantage that 
there is no trouble  even when we deal with singular varieties with $G$-action.

%%%%%%%%%%%%%%%%%%%%%%%%%%%%%%%%%%%%%%%%%%%
\section{Thom polynomials}
In this section we discuss on $G$-characteristic classes 
associated to $G$-classifications of a $G$-space 
(=classifications of ``singularities" of various objects).  
This correspondence is given by  a ``Segre-version" of $C^G_*$.  
We assume $k=\C$ in Subsection \ref{tp}, but in the other subsection 
we think of  both contexts although we write $H_*$ throughout. 
\subsection{Thom polynomial} \label{tp}
Let  $X$ be a nonsingular $G$-variety. 
Then, we have 
$${\rm Dual}_G\circ  C^G_*: \F^G_{inv}(X) \to H_G^*(X)$$
 by using $G$-Poincar\'e dual. 
For a $G$-invariant  subvariety $W$ of $X$ with codimension $l$, 
the leading term of ${\rm Dual}_G\circ C^{G}_*(\jeden_W)$  
is the $G$-equivariant Poincar\'e dual to $[W]_G$ in $X$, 
which is usually called {\it the Thom polynomial} of $W$ (in
$G$-classification of $X$), 
denoted by $tp(W) \in H_G^{2l}(X)$, 
cf.  \cite{Kaz1}, \cite{FR} and their references. 

In particular,  if $X$ is a $G$-affine space (as a usual case in $tp$ theory) , 
$H_G^*(X)=H^*(BG)$ and hence 
$tp(W)$ is written as a  polynomial of characteristic classes $c_i$ of
$G$-bundles, 
which has the  ``universality" in the following sense: 

\noindent
(Universality): For any  bundle  $E \to M$ with fibre $X$ and the structure
group $G$ 
over a nonsingular base space $M$ of dimension $m$, 
we associate a subbundle $E_W \to M$ with fibre $W$ (because $W$ is
$G$-invariant). 
For a ``generic" section $s:M \to E$, we set 
$$W(s):=s^{-1}(E_W)$$
and  call this  {\it the singular set of type $W$}, 
which has the expected codimension $l=\codim W$. 
Let $i:W(s) \to M$ be the inclusion. 
Then, the fundamental class of the singular set is expressed in $M$  by 
$$i_*[W(s)]=tp(W)(c(E))\frown [M] \in H_{2(m-l)}(M)$$
after substituting $c_i(E)$ to $c_i$ arising in $tp(W) \in H^{2l}(BG)$. 

This theorem  basically  goes back to R. Thom \cite{Thom}, 
and it is proved in a topological setting (cf. \cite{Kaz1}, \cite{FR}). 
In an algebraic setting, we show a more general statement  later. 

As a typical example of Thom polynomials, there is 
the so-called {\it Thom-Porteous formula}:  
Let $X$ be the affine space $Hom(\C^m, \C^{m+k})$  
on which the group $G=GL(m,\C)\times GL(n+k,\C)$ 
operates from the right and left as linear coordinate changes. 
The invariant subvariety $W$ under consideration is 
 the closure of the orbit  with kernel dimension $i$, usually denoted by
$\bar\Sigma^i$. 
Let $f: E\to F$  be a suitably generic vector bundle map over $M$, i.e., 
a section $f: M \to Hom(E,F)$,  
where $E$ and $F$ are of rank $m$ and $m+k$, respectively. 
Then the fundamental class of the degenerate loci $\bar\Sigma^i(f)$ 
 is expressed (as in cohomology of $M$)  
by a certain Schur polynomial in $c_i(F-E)$, that is 
 $tp(\bar\Sigma^i)$, cf.   \cite{F} Chap. 14.   
In this section concerning $tp$, 
 it would be helpful  to take this example in mind throughout.

\subsection{Generic morphisms with respect to a subvariety} \label{generic}  
Throughout this subsection, we forget the group action. 
Let $k=\C$ and $f$  be 
a morphism of a possibly singular variety  into a nonsingular variety,   
and $W$ a subvariety of the target variety of $f$. 
Note that any  corresponding complex analytic variety admits a Whitney
stratification. 
We say $f$ is transverse to  $W$  (or say, generic with respect to $W$) 
if the restriction of $f$ to any stratum of the source variety 
is transverse to any strata of $W$. 

Instead, we define the ``genericity" of $f$ with respect to $W$  in the
algebraic context.  
It is, in fact, related to the VRR theorem. 
Let us remind that 
the (ordinary) VRR theorem for {\it smooth morphisms} $f$  (Yokura \cite{Y1}) 
says that $C_*$ is compatible with $f^*$ and 
$f^{**}=c(\nu_f)^{-1} \frown f^*$, 
however the theorem fails 
for non-smooth morphisms (especially, regular embeddings).    
On one hand, 
morphisms $f$ in which we are now interested are 
regular embeddings, or  more suitably,  
local complete intersection morphisms,  
cf.  \cite{F}  
(a {\it l.c.i. morphisms} $f:X \to Y$ is a composition of a regular embedding  
$i: X \to N$ and a smooth morphism (e.g., fibrations)  $p: N \to Y$; 
it has the virtual normal bundle $\nu_f:=\nu_X-T_p$). 
Typical examples we take in mind 
are sections of  vector bundles (e.g., sections $X \to Y:=Hom(E,F)$) or 
any morphisms between nonsingular varieties.  

Although the VRR theorem for l.c.i morphisms fails, 
the VRR formula for a {\it distinguished} element $\jeden_W \in \F(Y)$ 
makes sense, that is our definition of 
 the ``genericity of $f$ with respect to $W$": 
\begin{definition} 
Let $f: X \to Y$  be a  l.c.i. morphism 
with the virtual normal bundle $\nu$. 
We say that $f$ is {\it generic} with respect to a subvariety $W$  of  $Y$ 
 if  it holds that $C_*\circ f^*(\jeden_W)=f^{**}\circ C_*(\jeden_W)$, i.e., 
$i_*C^{SM}(f^{-1}(W))=c(\nu)^{-1}\frown f^* i_*C^{SM}(W)$ 
where $i_*$ are induced maps via inclusions. 
\end{definition}
\begin{prop} 
(1) Any smooth morphism $f: X \to Y$ is 
generic with respect to any subvariety in $Y$.  
(2) Any l.c.i. morphism $f: X \to Y$ 
is generic with respect to $Y$. 
\end{prop}
\proof 
This is straightforward from the definition. %\qed 

\begin{rem} 
It was J.~ Sch\"urmann who established 
 {\it the generalized VRR formula for l.c.i. morphisms}  (Theorem 0.1
\cite{Schurmann}), 
which describes the defect $C_*\circ f^*-f^{**}\circ C_*$ in terms of 
 {\it the generalized vanishing cycle functor}.  
In the following subsections, 
we may state  theorems by using his vanishing cycle functor, 
instead of assuming  the ``genericity".  
We also remark that in  the complex case,  
the ``transversality" (in the sense of the stratification theory) implies  
the ``genericity" in  the above sense
(cf.  Proposition 1.3  of \cite{PP}, \cite{O1}, Corollary 0.1 of 
\cite{Schurmann}).  
\end{rem}

\subsection{Schwartz-MacPherson Segre class} \label{SMSegre} 
For a closed subscheme $Z$ in a {\it nonsingular} variety $M$, 
we define 
$$s^{SM}(Z,M):=c(TM|_{Z})^{-1}\frown C^{SM}(Z) \; \in \; H_*(Z)$$ 
as an analogy  to 
the relation of the Segre covariance class and Fulton's canonical class
defined in \cite{F}, cf. Remark \ref{ChernF}. 
This ``Segre version of SM classes" has been introduced also in \cite{Aluffi}.  
We may  denote this class by $s^{SM}(\jeden_Z,M)$  
($Z$ being  a (reduced) subvariety). 
We also define $s^{M}(Z,M)$ by replacing $C^{SM}$ to $C^{M}$ 
(for a subvariety $Z$,  $s^{M}(\jeden_Z,M)=s^{M}(Eu_Z,M)$). 

Let us return to our equivariant setting. We give a generalization of $tp$
as follows: 

\begin{definition} 
For an invariant subvariety $W$ in a nonsingular $G$-variety $X$ 
(its  $G$-inclusion is denoted by $i:W \to X$), 
we define  {\it the universal Segre-SM class} 
$$s^{SM}_G(W, X) :=c^G(TX|_W)^{-1} \frown  C^G_*(\jeden_W)  \; \in \;
H^G_*(W), $$
or equivalently, 
$s^{SM}_G(W, X)= \varphi_U s^{SM}(W\times_G U, X\times_G U)$ 
where $\varphi_U$ is the limit map $H_{trunc}(W\times_G U) \to H^G_*(W)$. 
Its $G$-equivariant dual in $X$ is denoted by 
$$tp^{SM}(W) := {\Dual}_G\; i^G_*s^{SM}_G(W, X) \; \in \;  H_G^*(X).$$
\end{definition}

Note that $tp^{SM}(W)$  is  a formal power series 
$$tp^{SM}(W)=\sum_{i=0}^\infty  tp_i^{SM}(W) \in H_G^*(X)=\prod H^i_G(X).$$

\subsection{Universality for sections} 
Let $s$ be a section of a bundle $E\to M$  with fibre $X$ and structure
group $G$, 
$W$  an invariant subvariety of $X$.  
For short, we say  $s:M \to E$ is  {\it generic} (with respect to $W$) 
if the morphism $s$ is generic with respect to the associated subbundle
$E_W$ with fibre $W$.  
In the following theorem, 
 we assume that $M$ is a  {\it quasi-projective} variety. 

\begin{thm} \label{segretp} 
Let $X$ be a $G$-affine space,  $W$ an invariant subvariety of $X$ 
of codimension $l$. Then, \\
(1)  $tp_i^{SM}(W)=0$ ($i < l$) and $tp_l^{SM}(W)$ coincides with the Thom
polynomial $tp(W)$: 
$$tp^{SM}(W)=tp(W) + \mbox{higher terms}.$$
(2) (universality) For any generic bundle $E\to M$ and any generic section
$s$ w.r.t. $W$, 
we have 
$$i_*C^{SM}(W(s))=tp^{SM}(W)(c(E)) \frown C^{SM}(M) \; \in \; H_*(M),$$ 
where $i_*$ denotes the induced map by the inclusion. 
\end{thm}

\begin{rem} As a generalization of Thom-Porteous formula, 
Parusi\'nski-Pragacz \cite{PP} give a formula of Schwartz-MacPherson classes of 
degeneracy loci of bundle maps, 
that is one of our motivation to define our generalized tp as given above. 
In fact, they computed  $tp^{SM}$  of $\bar\Sigma^{i}$, see Theorem 2.1 in
\cite{PP}. 
\end{rem}

\proof 
(1) This is obvious, in fact  the top term of 
$i^G_*\, C^{SM}_G(W)$ is just $i^G_*[W]_G \in H^G_{2(n-l)}(X)$.  

\t
(2) The main point is  the following key lemma (Lemma 1.6 of Totaro
\cite{Totaro})  
on the existence of classifying maps of $G$-bundles over a 
quasi-projective variety: 

\begin{lem} \label{totaro} (\cite{Totaro}): 
For any algebraic  bundle $E\to M$ 
with fibre $X$ and structure group $G$ over a quasi-projective variety $M$, 
there is a bundle $q:M_1 \to M$ with fibre  being an affine space 
which admits an algebraic {\it classifying map} $\rho: M_1 \to U/G$ 
for some large $U=V-S$ ($\in I(G)$) 
so that $\rho^*(X\times_G U)\simeq q^*E$. 
$$
\begin{array}{ccccccc}
E & \longleftarrow & q^*E & 
\overset{\bar\rho}{\longrightarrow}& X\times_G U& \to & X\times_G EG\\
\downarrow && \downarrow & & \downarrow && \downarrow\\
M& \overset{q}{\longleftarrow} & M_1 &\overset{\rho}{\longrightarrow} & U/G
 &  \to & BG
\end{array}
$$
\end{lem}

The rest of the proof is straightforward from this lemma, the inductive
limit argument 
and the genericity of morphisms w.r.t. certain varieties associated to $W$. 

We take $M_1$ as in Lemma \ref{totaro} and denote $M_2:=X\times_G U$ and 
$W_G:=W\times_G U$ ($\subset M_2$). 
Now $X$ is assumed to be an affine space,  so $M_2$ is nonsingular, and hence 
$C^{SM}(M_2)=c(TM_2)\frown [M_2]$. 
Then 
\begin{eqnarray*}
&&r_U \, tp^{SM}(W) \frown C^{SM}(M_2) \\
&&\qquad \quad = r_U \, tp^{SM}(W) \cdot  c(TM_2) \frown [M_2] \\
&&\qquad \quad =c(TM_2)\cdot r_U \circ {\rm Dual}_G \circ  i^G_*
\left(s^{SM}_G(W,X)\right) \frown [M_2] \\
&&\qquad \quad =c(TM_2)\cdot  {\rm Dual}^{-1}\circ r_U\circ {\rm Dual}_G \circ
\varphi_U(i_*s^{SM}(W_G,M_2))\\
&&\qquad \quad =c(TM_2)\frown i_*s^{SM}(W_G,M_2)\\
&&\qquad \quad =c(TM_2)\frown \left(c(TM_2)^{-1} \frown i_*C^{SM}(W_G)\right)\\
&&\qquad \quad =  i_*\,C^{SM}(W_G). 
\end{eqnarray*}
The section $s:M \to E$ induces a section $s':M_1 \to q^*E$, and then set 
$f:=\bar\rho\circ s':M_1 \to M_2$ and $\nu_f$ the virtual normal bundle of $f$ 
(in fact $\bar\rho$ can be taken as a regular embedding). 
Since we assume $s$ is generic (w.r.t. $E_W$), 
it turns out that $f$ is generic w.r.t. $W_G$ and hence we have 
$$i_*\,C^{SM}(W(s'))=c(\nu_f)^{-1}\frown f^*\, i_*\,C^{SM}(W_G), $$
$$C^{SM}(M_1)=c(\nu_f)^{-1}\frown f^*\,  C^{SM}(M_2).$$
Also for the smooth morphism $q:M_1 \to M$, 
$$i_*\,C^{SM}(W(s'))= q^{**}\, i_*\,C^{SM}(W(s)), \quad 
C^{SM}(M_1)= q^{**}\, C^{SM}(M),$$
where $q^{**}=c(T_q)\frown q^*:H_*(M) \to H_*(M_1)$. 
Since $H_G^*(X)=H_G^*(pt) \simeq H^*(BG)$ ($X$ being a $G$-affine space), 
$f^*r_U$ is identified with the pullback via the classifying map $\rho^*r_U$, 
which sends 
the universal characteristic class $c_i \in H^{2i}(BG)$ 
to the Chern class $c_i(q^*E) \in H^{2i}(M_1)$.  
Thus, we have 
\begin{eqnarray*}
q^{**}\,i_*\,C^{SM}(W(s)) &=& i_*\,C^{SM}(W(s'))  \\
&=& c(\nu_f)^{-1}\frown f^*\, i_*\,C^{SM}(W_G) \\
&=& c(\nu_f)^{-1} \frown f^*\left(r_U \,tp^{SM}(W)\frown C^{SM}(M_2) \right)\\
&=& f^*\,r_U \,tp^{SM}(W)
\frown \left(c(\nu_f)^{-1}\frown f^*\,C^{SM}(M_2)\right)\\
&=&  tp^{SM}(W)(c(q^*E))\frown C^{SM}(M_1) \\
&=&  tp^{SM}(W)(c(q^*E))\frown q^{**}\, C^{SM}(M) \\
&=&  q^{**}\,\left(tp^{SM}(W)(c(E))\frown C^{SM}(M) \right). 
\end{eqnarray*}
Since $q^{**}$ is an isomorphism,  this equality shows (2).  
%\qed

\

\begin{cor} Let $k=\C$. For a generic section $s$ as above, 
the topological Euler characteristic of the singular set $W(s)$ of type $W$ 
is given by 
$$\chi(W(s))=\int_M tp^{SM}(W)(c(E)) \frown C^{SM}(M).$$ 
\end{cor}

\proof 
This is a direct consequence from the fact that 
the Euler characteristic of a variety  is equal to 
the $0$-th degree of its (ordinary) Schwartz-MacPherson class.  %\qed 

\begin{rem}
We summarize this section.  
As seen in \S 6, our equivariant MacPherson theory $C^G_*$ certainly 
gives the Chern class version of  
various Euler characteristics of quotients, 
while in this section, 
viewing from the top dimensional side, 
our $C^G_*$ (precisely its Segre-version)  gives 
``higher dimensional" Thom polynomial theory.  

For a nonsingular $G$-variety $X$, 
we have introduced universal Segre-SM classes of invariant subvarieties, 
that produces a homomorphism between abelian groups 
$$tp^{SM}: \F^G_{inv}(X) \to H^*_G(X).$$ 
This is actually a natural transformation 
for the category of nonsingular $G$-varieties and proper $G$-morphisms. 
In particular,  if $X$ be a $G$-affine space, 
the values are universal polynomials in Chern classes for $G$-bundles. 
In other words, $tp^{SM}$ is a correspondence 
from a ``{\it local $G$-classification}"  of singularities 
to ``{\it global invariants}" for 
generic sections of any associated affine bundles $E \to M$, 
that we called   the``{\it universality}" of $tp^{SM}$ throughout this section. 
The property is shown totally 
in the algebraic context, by using Totaro's classifying maps.  
Besides, we may define $tp^{M}$ and $tp^F$ in the same way as $tp^{SM}$, 
the relations among which should be related to 
local invariants of closures of $G$-orbits (e.g., local Euler obstruction). 

As a typical example, let us explain this ``local-global" correspondence in 
Singularity theory of differentiable mappings between manifolds: 
The classification problem of map-germs is reduced to 
the classification in the level of some $r$-jets, i.e., the so-called {\it finite determinacy}. 
So then we think of  some $r$-jet space $X=J^r(m, m+k)$ 
of germs $\C^m,0 \to \C^{m+k}, 0$ 
(the space of $r$-th Taylor expansions) with the action 
of the group $G$ of $r$-jets of coordinate changes of source and target 
(called the $\cal RL$ (right-left) classification).  
An invariant subvariety (e.g., the closure of an orbit or of a family of orbits) 
 is called {\it a singularity type}.  
The case of $r=1$ is just singularities $\bar{\Sigma}^i$ of vector bundle maps. 
For any suitable generic map $f:M^m \to N^{m+k}$, 
a singularity type $\eta$ in the local classification gives a ``global invariant of $f$" 
$tp(\eta)$ and $tp^{SM}(\eta)$,  which universally express 
the fundamental class and the Chern class of $\eta(f)$ 
(Precisely saying, the {\it $\K$-classification} is much useful 
for the classification of generic singularities of maps, 
then  a further inductive limit process arises 
according to the parameters $r \to \infty$ and $m \to \infty$ 
with fixed difference $k$ of the dimension of source and target spaces. 
It turns out that $tp$ (each component of $tp^{SM}$) 
is a polynomial in $c_i$'s which mean Chern classes of the virtual normal bundle 
$c_i(f^*TN-TM)$ for any generic $f:M \to N$,  
see \cite{Kaz1}, \cite{FR}, \cite{O1}).  
There may be many interesting directions for further researches, 
for instance, 
computational aspects of $tp^{SM}$, $tp^M$ and $tp^F$  in relation 
with local invariants such as Milnor number and Polar multiplicities, 
and the $tp^{SM}$-version for multi-singularities (\cite{Kaz2})  
combined with materials in \S 6, etc. 
\end{rem}

%%%%%%%%%%%%%%%%%%%%%%%%%%%%%%%%%%%%%%%%%%%

\end{document}